\documentclass[12pt]{amsart}
\usepackage{latexsym}

\newcommand\ff{\varphi}
\newcommand\cA{\mathcal{A}}
\newcommand\CC{\mathbf{C}}
\newcommand\NN{\mathbf{N}}
\newcommand\la{\langle}
\newcommand\ra{\rangle}

\newcommand\entspricht{\hat=}

\newfont{\bs}{cmbxsl10 scaled 1200}

\begin{document}

\title{Combinatorics of free cumulants}


\author{Bernadette Krawczyk}
\author{Roland Speicher}
\email{roland.speicher@urz.uni-heidelberg.de}
\thanks{The second author was supported by a Heisenberg fellowship of the
DFG}
\address{Institut f\"ur Angewandte Mathematik, Im Neuenheimer Feld
294, D-69120 Heidelberg, Germany}
\begin{abstract}
We derive a formula for expressing free cumulants whose
entries are products of random variables in terms of the lattice
structure of non-crossing partitions.
We show the usefulness of that result by giving direct and conceptually
simple proofs for a lot of results about $R$-diagonal elements.
Our
investigations do not assume the trace property for the considered linear
functionals.
\end{abstract}
\maketitle

\section*{Introduction}
Free probability theory, due to Voiculescu \cite{V,VDN},
is a non-com\-muta\-tive probability theory where
the classical concept of ``independence" is replaced by a
non-com\-muta\-tive
analogue, called ``freeness". Originally this theory was introduced
in an operator-algebraic context for dealing with questions on special
von Neumann algebras. However, since these beginnings free probability
theory has evolved into a theory
with a lot of links to quite different fields. In particular, there
exists a combinatorial facet: main aspects of free probability theory
can be considered as the combinatorics of non-crossing partitions.

There are two main approaches to freeness:
\begin{itemize}
\item
the original approach, due
to Voiculescu, is analytical in nature and relies on
special Fock space constructions for the considered distributions.
\item
the approach of Speicher \cite{S2,S1,S3} is combinatorial in nature and
describes freeness in terms of so-called free cumulants -- these objects
are defined via a precise combinatorial description involving the lattice
of non-crossing partitions; a lot of questions on freeness reduce in
this approach finally to combinatorial problems on non-crossing partitions.
\end{itemize}
The relation between these two approaches is given by the fact that
the free cumulants appear as coefficients in the operators constructed
in the Fock space approach. This connection was worked out by Nica
\cite{N}.

Here, we will investigate one fundamental problem in the combinatorial
approach and show that there is a beautiful combinatorial structure behind
this.

In the combinatorial approach to freeness one defines,
for a given linear functional $\ff$ on a unital algebra $\cA$,
so-called free cumulants
$k_n$ ($n\in\NN$), where each $k_n$ is a multi-linear functional
on $\cA$ in $n$
arguments. The connection between $\ff$ and the $k_n$ is given by
a combinatorial formula involving the lattice of non-crossing partitions.
(The name ``cumulants" comes from classical probability
theory; there exist analogous objects with that name, the only difference
is that there all partitions instead of non-crossing partitions appear.)
It seems that many problems on freeness are easier to handle in terms
of these free cumulants than in terms of moments of $\ff$.
In particular, the definition of freeness itself becomes much handier for
cumulants than for moments.
Since cumulants are multi-linear objects this implies that
for problems involving the linear structure of the algebra $\cA$ cumulants
are quite easily and effectively to use.
For problems involving the multiplicative structure of $\cA$, however,
it is not so clear from the beginning that cumulants are a useful tool
for such investigations. Nevertheless in a lot of examples
it has turned out that this is indeed the case.
In a sense, we will here present the unifying reason for these positive
results.
Namely, dealing with multiplicative problems reduces on the level of
cumulants
essentially to the problem of understanding the structure
of cumulants whose arguments are products of variables. Here, in Section 2,
we will
show that this can be understood quite well and that there exists a
nice and simple combinatorial description for such cumulants.

That this formula is also useful will be demonstrated in Section 3. We
will reprove and generalize a lot of results around the multiplication
of free random variables. In particular, we will consider an
important special class
of distributions, so-called $R$-diagonal elements. These were introduced
by Nica and Speicher in \cite{NS1}. However,
the investigations and characterizations in \cite{NS1,NS3} were not always
straightforward and used a lot of ad hoc combinatorics.
Our approach here will be much more direct and conceptually
clearer. Furthermore, we will get in the same spirit direct proofs of
results of Haagerup and Larsen \cite {HL,L}
on powers of $R$-diagonal elements.

An important point to make is that all earlier investigations
on $R$-diagonal elements were
always restricted to a tracial frame -- i.e., $\ff$ was assumed to satisfy
the trace condition
$\ff(ab)=\ff(ba)$ for all $a,b\in\cA$.
In contrast, our approach does not rely on this assumption, so all our
results are also valid for non-tracial $\ff$. Thus we do not only
get simple proofs for known results but also generalizations of all
these results to the general, non-tracial case. (That
non-tracial
$R$-diagonal elements appear quite naturally
can, e.g., be seen in \cite{S},
where such elements arise in the polar decomposition of
generalized circular elements).

Our Propositions 3.5 and 3.9 were inspired by and prove some conjectures of
the recent work \cite{NSS}.
There the notion of $R$-diagonality is also treated in the non-tracial case
and some of our results
of Section 3 are proved there for the
general case, too. However, the approach in \cite{NSS} is quite different
from
the present one
and relies on Fock space representations and freeness with amalgamation.

The paper is organized as follows. In Section 1, we give a short and
self-contained summary of the relevant basic definitions and facts about
free
probability theory and non-crossing partitions. In Section 2, we state
and prove our main combinatorial result on the structure of free cumulants
whose arguments are products and, in Section 3, we apply this result to
derive various statements about $R$-diagonal elements.

\section{Preliminaries}

In this section we provide a short and self-contained summary of
the basic definitions and facts needed for our later investigations.

\subsection{Non-commutative probability theory.}
1) We will always work in the frame of a
{\it non-commutative probability space} $(\mathcal{ A},\varphi)$.
This is, by definition, a
pair consisting of a unital $\ast$-algebra $\mathcal{ A}$ and a
unital linear
functional $ \varphi \: : \mathcal{ A} \to \mathbf{C}$.
($\varphi$ unital means that
$\varphi(1) =1$.)
\\
The elements $a \in \mathcal{ A}$ are called {\it non-commutative
random variables}, or just {\it random variables}
in $(\mathcal{ A},\varphi)$.
\\
Let $a_1,
\dots, a_n$ be random variables in a non-com\-mu\-ta\-tive
probability space $(\mathcal{ A},\varphi)$.
Let $\CC\la X_1,\dots,X_n\ra$ denote the algebra of polynomials in
$n$ non-commuting indeterminants -- i.e., the algebra generated by
$n$ free generators. Then the linear functional
$$\mu_{a_1, \dots, a_n}:\CC\la X_1,\dots,X_n\ra\to\CC$$
given by linear extension of
$$X_{i(1)}\dots X_{i(m)}\mapsto \varphi(a_{i(1)} \cdots a_{i(m)})
\qquad(m \in {\mathbf{N}}, \: 1 \leq i(1),\dots,i(m)\leq n)
$$
is called the {\it joint distribution} of $a_1, \dots, a_n$.
\\
The joint distribution of $a$ and $a^*$ is also called
the {\it $*$-distribution} of $a$.
\\
Consider random variables $a_i$ and
$b_i \: (1 \leq i \leq n)$ in $(\mathcal{ A},\varphi)$. Then
$a_1,\dots,a_n$ and $b_1,\dots,b_n$ have the same joint
distribution, if the following equation holds for all $m \in
{\mathbf{N}}, \: 1 \leq i(1),\dots,i(m) \leq n$:
\[ \varphi(a_{i(1)} \cdots a_{i(m)})
   = \varphi(b_{i(1)} \cdots b_{i(m)}) \: .
\]
2) Note that all our considerations will be on the algebraic (or
combinatorial) level, thus we will not require that $\varphi$ is
a positive functional.
However, it is well known that freeness --
the crucial structure in our investigations - is compatible with
positivity properties. The requirement that our probability space
should be a $*$-algebra and not just an algebra is only for convenience,
since, in Section 3, we will need the $*$ for dealing with Haar
unitaries and $R$-diagonal elements. In all statements where no $*$
appears we could also replace the requirement ``$*$-algebra" by
``algebra".
\\
3) Most of the questions which we will investigate in Section 3 were
up to now only considered for tracial linear functionals. We stress that
all our considerations do not use the trace property, i.\/e.\/
we will not use the equation
$\varphi(ab) =\varphi(ba)$.

\subsection{Partitions.}
1) Fix $n\in\NN$.
We call $\pi = \{ V_1, \dots, V_r \}$ a {\it partition}
of $S=(1,\dots,n)$ if and only if the $V_i \: (1 \leq i \leq r)$ are
pairwisely disjoint, non-void tuples such that $V_1 \cup \dots
\cup V_r = S$.
We call the tuples $V_1, \dots, V_r$ the {\it blocks} of
$\pi$. The number of components of a block $V$
is denoted by $|V|$. Given two
elements $p$ und $q$ with $1 \leq p,q \leq n$, we write $p
\sim_{\pi} q$, if $p$ and $q$ belong to the same block of $\pi$.
\\
We get a linear representation of a partition $\pi$ by writing
all elements $1, \dots, n$ in a line, supplying each with a
vertical line under it and joining the vertical lines of the
elements in the same block with a horizontal line.
\\
Example: A partition of the tuple $S =
(1,2,3,4,5,6,7)$ is
\\
\vskip0.0cm
\setlength{\unitlength}{0.3cm}
\[ \pi_1 = \{ (1,4,5,7), (2,3), (6) \}
   \qquad \entspricht \qquad
   \begin{picture}(6,2)\thicklines
   \put(0,0){\line(0,1){2}}
   \put(0,0){\line(1,0){6}}
   \put(3,0){\line(0,1){2}}
   \put(4,0){\line(0,1){2}}
   \put(6,0){\line(0,1){2}}
   \put(1,1){\line(0,1){1}}
   \put(2,1){\line(0,1){1}}
   \put(1,1){\line(1,0){1}}
   \put(5,1){\line(0,1){1}}
   \put(-0.3,2.2){1}
   \put(0.7,2.2){2}
   \put(1.7,2.2){3}
   \put(2.7,2.2){4}
   \put(3.7,2.2){5}
   \put(4.7,2.2){6}
   \put(5.7,2.2){7}
   \end{picture} \: .
\]
If we write a block $V$ of a partition in the form $V=(v_1,\dots,v_p)$
then this shall always imply that $v_1<v_2<\dots<v_p$.
\\
2) A partition $\pi$ is called {\it non-crossing}, if the following
situation does not occur: There exist $1 \leq p_1 < q_1 < p_2 <
q_2 \leq n$ such that ${p_1} \sim_{\pi} {p_2} \not\sim_{\pi}
{q_1} \sim_{\pi} {q_2}$:
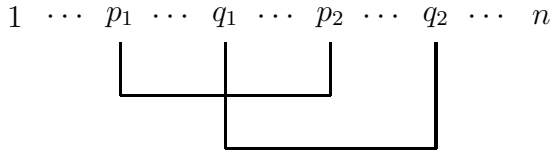
\setlength{\unitlength}{0.7cm}
\[ \begin{picture}(10,3)\thicklines
   \put(2,1){\line(0,1){1}}
   \put(2,1){\line(1,0){4}}
   \put(6,1){\line(0,1){1}}
   \put(4,0){\line(0,1){2}}
   \put(4,0){\line(1,0){4}}
   \put(8,0){\line(0,1){2}}
   \put(0,2.5){\makebox(0,0){$1$}}
   \put(1,2.5){\makebox(0,0){$\cdots$}}
   \put(2,2.5){\makebox(0,0){${p_1}$}}
   \put(3,2.5){\makebox(0,0){$\cdots$}}
   \put(4,2.5){\makebox(0,0){${q_1}$}}
   \put(5,2.5){\makebox(0,0){$\cdots$}}
   \put(6,2.5){\makebox(0,0){${p_2}$}}
   \put(7,2.5){\makebox(0,0){$\cdots$}}
   \put(8,2.5){\makebox(0,0){${q_2}$}}
   \put(9,2.5){\makebox(0,0){$\cdots$}}
   \put(10,2.5){\makebox(0,0){$n$}}
   \end{picture}
\]
The set of all non-crossing partitions of $(1,\dots,n)$ is
denoted by $NC(n)$.
In the same way as for $(1,\dots,n)$ one can introduce non-crossing
partitions $NC(S)$ for each finite linearly ordered set $S$. Of course,
$NC(S)$ depends only on the number of elements in $S$. In our
investigations,
non-crossing partitions will appear as partitions of the index set of
products of random variables $a_1\cdots a_n$. In such a case, we will also
sometimes use the notation $NC(a_1,\dots,a_n)$.
(If some of the $a_i$ are equal, this might make no rigorous sense,
but there should arise no problems by this.)
\\
If $S$ is the union of two disjoint sets $S_1$ and $S_2$ then, for
$\pi_1\in NC(S_1)$ and $\pi_2\in NC(S_2)$, we let $\pi_1\cup \pi_2$
be that partition of $S$ which has as blocks the blocks of $\pi_1$ and
the blocks of $\pi_2$. Note that $\pi_1\cup\pi_2$ is not automatically
non-crossing.
\\
3) Let $\pi, \sigma \in NC(n)$ be two non-crossing partitions. We
write $\sigma \leq \pi$, if every block of $\sigma$ is completely
included in a block of $\pi$. Hence, we obtain $\sigma$ out of
$\pi$ by refining the block-structure.
For example, we have
$$\{(1,3),(2),(4,5),(6,8),(7)\}\leq \{(1,3,7),(2),(4,5,6,8)\}.$$
The partial order $\leq$ induces a lattice structure on
$NC(n)$.
In particular, given two non-crossing partitions
$\pi,\sigma\in NC(n)$,
we have their join
$\pi \vee \sigma$, which is the unique
smallest $\tau \in NC(n)$ such that $\tau \geq \pi$
and $\tau \geq \sigma$.
\\
The maximum of
$NC(n)$ -- the partition which consists of one block with $n$
components -- is denoted by $1_n$. The partition consisting of
$n$ blocks, each of which has one component, is the minimum of
$NC(n)$ and denoted by $0_n$.
\\
4) The lattice $NC(n)$ is self-dual and there exists an
important anti-isomorphism $K:NC(n)\to NC(n)$ implementing this
self-duality. This complementation map $K$ is defined as follows:
Let $\pi$ be a non-crossing partition of the numbers
$1,\dots,n$. Furthermore, we consider numbers
$\bar{1},\dots,\bar{n}$ with all numbers ordered like
\[ 1 \: \bar{1} \: 2 \: \bar{2} \dots n \: \bar{n} \: .
\]
The {\it complement} $K(\pi)$ of $\pi \in NC(n)$ is defined to be the
biggest $\sigma \in NC(\bar{1},\dots,\bar{n})\entspricht NC(n)$ with
\[ \pi \cup \sigma \in NC(1,\bar{1},\dots,n,\bar{n}) \: .
\]
Example: Consider the partition
$\pi := \{ (1,2,7), (3), (4,6), (5), (8) \}\in NC(8)$. For the complement
$K(\pi)$ we get
\[ K(\pi) = \{ ({1}), ({2}, {3}, {6}), ({4},
                 {5}), ({7}, {8}) \} \: ,
\]
as can be seen from the graphical representation:
\\
\vskip5pt
\begin{center}
  \setlength{\unitlength}{0.3cm}
  \begin{picture}(16,6)
  \thicklines
  \put(1,0){\line(0,1){5}}
  \put(1,0){\line(1,0){12}}
  \put(3,0){\line(0,1){5}}
  \put(13,0){\line(0,1){5}}
  \put(15,1){\line(0,1){4}}
  \put(5,2){\line(0,1){3}}
  \put(7,2){\line(0,1){3}}
  \put(7,2){\line(1,0){4}}
  \put(11,2){\line(0,1){3}}
  \put(9,4){\line(0,1){1}}
  \linethickness{0.6mm}
  \put(2,1){\line(0,1){4}}
  \put(4,1){\line(1,0){8}}
  \put(4,1){\line(0,1){4}}
  \put(6,1){\line(0,1){4}}
  \put(12,1){\line(0,1){4}}
  \put(8,3){\line(0,1){2}}
  \put(8,3){\line(1,0){2}}
  \put(10,3){\line(0,1){2}}
  \put(14,0){\line(0,1){5}}
  \put(14,0){\line(1,0){2}}
  \put(16,0){\line(0,1){5}}
  \put(0.8,5.5){1}
  \put(1.8,5.5){$\bar 1$}
  \put(2.8,5.5){2}
  \put(3.8,5.5){$\bar 2$}
  \put(4.8,5.5){3}
  \put(5.8,5.5){$\bar 3$}
  \put(6.8,5.5){4}
  \put(7.8,5.5){$\bar 4$}
  \put(8.8,5.5){5}
  \put(9.8,5.5){$\bar 5$}
  \put(10.8,5.5){6}
  \put(11.8,5.5){$\bar 6$}
  \put(12.8,5.5){7}
  \put(13.8,5.5){$\bar 7$}
  \put(14.8,5.5){8}
  \put(15.8,5.5){$\bar 8$}
  \end{picture}
\hspace{0.5cm} .
\end{center}
5) Non-crossing partitions and the complementation map
were introduced by Kreweras \cite{K}; for
further combinatorial investigations on that lattice, see, e.g.,
\cite{E,SU}.
\\
6) The main combinatorial ingredient of Theorem 2.2 will be joins
with special partitions $\sigma$ whose blocks consist of neighbouring
elements, like $\pi \vee
\{(1),(2),\dots,(l, \dots, {l+k}),\dots,(n)\}$. This is
given by uniting the blocks of $\pi$
containing the elements $l, \dots,
{l+k}$, and we say that we obtain
$\pi \vee
\{(1),(2),\dots,(l, \dots, {l+k}),\dots,(n)\}$
by {\it connecting} the elements $l, \dots,
{l+k}$.
\\
Example: Considering the partition
\setlength{\unitlength}{0.3cm}
\[ \pi=\{ (1,8),(2,3),(4,5,7), (6) \}
   \qquad \entspricht \qquad
   \begin{picture}(7,3)\thicklines
   \put(0,0){\line(0,1){3}}
   \put(0,0){\line(1,0){7}}
   \put(7,0){\line(0,1){3}}
   \put(1,1){\line(0,1){2}}
   \put(1,1){\line(1,0){1}}
   \put(2,1){\line(0,1){2}}
   \put(3,1){\line(0,1){2}}
   \put(3,1){\line(1,0){3}}
   \put(4,1){\line(0,1){2}}
   \put(6,1){\line(0,1){2}}
   \put(5,2){\line(0,1){1}}
   \put(-0.3,3.2){1}
   \put(0.7,3.2){2}
   \put(1.7,3.2){3}
   \put(2.7,3.2){4}
   \put(3.7,3.2){5}
   \put(4.7,3.2){6}
   \put(5.7,3.2){7}
   \put(6.7,3.2){8}
   \end{picture}
\]
we have
\begin{align*}
\pi \vee \{(1,2,3,4),(5),(6),(7),(8)\} &= \{ (1,2,3,4,5,7,8), (6) \}\\
& \ \\
&\entspricht \qquad
   \begin{picture}(7,2)\thicklines
   \put(0,0){\line(0,1){2}}
   \put(0,0){\line(1,0){7}}
   \put(1,0){\line(0,1){2}}
   \put(2,0){\line(0,1){2}}
   \put(3,0){\line(0,1){2}}
   \put(4,0){\line(0,1){2}}
   \put(6,0){\line(0,1){2}}
   \put(7,0){\line(0,1){2}}
   \put(5,1){\line(0,1){1}}
   \put(-0.3,2.2){1}
   \put(0.7,2.2){2}
   \put(1.7,2.2){3}
   \put(2.7,2.2){4}
   \put(3.7,2.2){5}
   \put(4.7,2.2){6}
   \put(5.7,2.2){7}
   \put(6.7,2.2){8}
   \end{picture}
   \: .
\end{align*}

\subsection{Free cumulants.} \label{1.3}
Given a unital linear functional $\varphi:\cA\to\CC$ we define
corresponding
{\it (free) cumulants} $(k_n)_{n\in\NN}$
\begin{align*}
 k_n:\qquad\qquad \mathcal{ A}^n &\to {\mathbf{C}}, \\
      (a_1,\dots,a_n) &\mapsto k_n(a_1,\dots,a_n)
\end{align*}
indirectly by the following system of equations:
\begin{equation} \label{1.2,1}
   \varphi(a_1 \cdots a_n)
   = \sum_{\scriptstyle \pi \in NC(n)}
     k_{\pi} \left[a_1,\dots,a_n\right]
     \qquad (a_1,\dots,a_n \in \mathcal{ A}) \: ,
\end{equation}
where $k_\pi$ splits multiplicatively in a product of cumulants
according to the block structure of $\pi$, i.e.
\begin{equation} \label{1.2,2}
   k_{\pi} [ a_1,\dots,a_n ]
   := \prod_{i=1}^r k_{|V_i|} (a_{i,1}, \dots,a_{i,\vert V_i\vert})
\end{equation}
for a partition $\pi = \{ V_1, \dots, V_r \} \in NC(n)$
consisting of $r$ blocks of the form $V_i = (a_{i,1}, \dots,
a_{i,\vert V_i\vert})$.\\
The defining relation (1) expresses the moment
$\ff(a_1\cdots a_n)$ in terms of cumulants, but by induction this
can also be resolved for giving the cumulants uniquely in terms of
moments:
\begin{equation} \label{1.2,3}
   k_n(a_1,\dots,a_n)
   = \varphi(a_1 \cdots a_n)
      - \sum_{\pi \in NC(n) \atop \pi \neq 1_n}
     k_{\pi} \left[a_1,\dots,a_n\right]
\end{equation}
Since, by induction, we know all cumulants of smaller order, i.e., all
$k_{\pi} [a_1,\dots,a_n]$ for $\pi \in
NC(n)$ with $\pi \neq 1_n$, this leads to an expression for $k_n$
in terms of moments. Abstractly, this is, of course, just the Moebius
inversion of relation (\ref{1.2,1}) and has the following form
\begin{equation}
k_n(a_1,\dots,a_n)=\sum_{\pi\in NC(n)} \mu(\pi,1_n) \ff_\pi[a_1,\dots,
a_n],
\end{equation}
where $\mu$ is the Moebius function of the lattice of non-crossing
partitions and
where $\ff_\pi$ is defined in the same multiplicative way as $k_\pi$ if we
put
$\ff_n(a_1,\dots,a_n):=\ff(a_1\cdots a_n).$
\\
Examples: Let us give the concrete form of
$k_n(a_1,\dots,a_n)$ for $n=1,2,3$.
\begin{itemize}
\item $n=1$:
\setlength{\unitlength}{0.2cm}
  \newsavebox{\NCianew}
  \savebox{\NCianew}(0,1){
  \thicklines
  \put(0,0){\line(0,1){1}}}

\[ k_1(a_1) = \varphi(a_1) \: .
\]

\setlength{\unitlength}{0.2cm}
  \newsavebox{\NCiianew}
  \savebox{\NCiianew}(1,1){
  \thicklines
  \put(0,0){\line(0,1){1}}
  \put(0,0){\line(1,0){1}}
  \put(1,0){\line(0,1){1}}}

  \newsavebox{\NCiibnew}
  \savebox{\NCiibnew}(1,1){
  \thicklines
  \put(0,0){\line(0,1){1}}
  \put(1,0){\line(0,1){1}}}

\setlength{\unitlength}{0.3cm}
  \newsavebox{\NCiibnewa}
  \savebox{\NCiibnewa}(1,1){
  \thicklines
  \put(0,0){\line(0,1){1}}
  \put(1,0){\line(0,1){1}}}

\item $n=2$: The only partition $\pi \in NC(2), \pi \neq 1_2$ is
$\usebox{\NCiibnewa} \quad .$ So we get
\begin{eqnarray*}
k_2(a_1,a_2)
& = & \varphi(a_1a_2) - k_{\usebox{\NCiibnew}\,\,\,}[a_1,a_2] \\
& = & \varphi(a_1a_2) - k_1(a_1)k_1(a_2) \\
& = & \varphi(a_1a_2) - \varphi(a_1)\varphi(a_2).
\end{eqnarray*}
Using the notation $\ff_\pi$ we can
also write this as
$$k_2(a_1,a_2)=
 \varphi_{\usebox{\NCiianew}\,\,\,}[a_1,a_2]
      - \varphi_{\usebox{\NCiibnew}\,\,\,}[a_1,a_2].$$

\item $n=3$: We have to take all partitions in $NC(3)$ except $1_3$,
i.\/e., the following partitions:

\setlength{\unitlength}{0.2cm}
  \newsavebox{\NCiiiannewa}
  \savebox{\NCiiiannewa}(2,2){
  \thicklines
  \put(0,0){\line(0,1){2}}
  \put(0,0){\line(1,0){2}}
  \put(1,0){\line(0,1){2}}
  \put(2,0){\line(0,1){2}}}

  \newsavebox{\NCiiibnnewa}
  \savebox{\NCiiibnnewa}(2,2){
  \thicklines
  \put(0,0){\line(0,1){2}}
  \put(1,0){\line(0,1){2}}
  \put(1,0){\line(1,0){1}}
  \put(2,0){\line(0,1){2}}}

  \newsavebox{\NCiiicnnewa}
  \savebox{\NCiiicnnewa}(2,2){
  \thicklines
  \put(0,0){\line(0,1){2}}
  \put(0,0){\line(1,0){1}}
  \put(1,0){\line(0,1){2}}
  \put(2,0){\line(0,1){2}}}

  \newsavebox{\NCiiidnnewa}
  \savebox{\NCiiidnnewa}(2,2){
  \thicklines
  \put(0,0){\line(0,1){2}}
  \put(0,0){\line(1,0){2}}
  \put(2,0){\line(0,1){2}}
  \put(1,1){\line(0,1){1}}}

  \newsavebox{\NCiiiennewa}
  \savebox{\NCiiiennewa}(2,2){
  \thicklines
  \put(0,0){\line(0,1){2}}
  \put(1,0){\line(0,1){2}}
  \put(2,0){\line(0,1){2}}}

\[ \usebox{\NCiiibnnewa} \quad ,
   \qquad \usebox{\NCiiicnnewa} \quad ,
   \qquad \usebox{\NCiiidnnewa} \quad ,
   \qquad \usebox{\NCiiiennewa} \quad .
\]

\setlength{\unitlength}{0.1cm}
  \newsavebox{\NCiiiannew}
  \savebox{\NCiiiannew}(2,2){
  \thicklines
  \put(0,0){\line(0,1){2}}
  \put(0,0){\line(1,0){2}}
  \put(1,0){\line(0,1){2}}
  \put(2,0){\line(0,1){2}}}

  \newsavebox{\NCiiibnnew}
  \savebox{\NCiiibnnew}(2,2){
  \thicklines
  \put(0,0){\line(0,1){2}}
  \put(1,0){\line(0,1){2}}
  \put(1,0){\line(1,0){1}}
  \put(2,0){\line(0,1){2}}}

  \newsavebox{\NCiiicnnew}
  \savebox{\NCiiicnnew}(2,2){
  \thicklines
  \put(0,0){\line(0,1){2}}
  \put(0,0){\line(1,0){1}}
  \put(1,0){\line(0,1){2}}
  \put(2,0){\line(0,1){2}}}

  \newsavebox{\NCiiidnnew}
  \savebox{\NCiiidnnew}(2,2){
  \thicklines
  \put(0,0){\line(0,1){2}}
  \put(0,0){\line(1,0){2}}
  \put(2,0){\line(0,1){2}}
  \put(1,1){\line(0,1){1}}}

  \newsavebox{\NCiiiennew}
  \savebox{\NCiiiennew}(2,2){
  \thicklines
  \put(0,0){\line(0,1){2}}
  \put(1,0){\line(0,1){2}}
  \put(2,0){\line(0,1){2}}}

With this we obtain:
\begin{eqnarray*}
k_3(a_1,a_2,a_3)
& = & \varphi(a_1a_2a_3)
      - k_{\usebox{\NCiiibnnew} \,\,\,}[a_1,a_2,a_3]
      - k_{\usebox{\NCiiicnnew} \,\,\,}[a_1,a_2,a_3] \\
  & & - k_{\usebox{\NCiiidnnew} \,\,\,}[a_1,a_2,a_3]
      - k_{\usebox{\NCiiiennew} \,\,\,}[a_1,a_2,a_3] \\
& = & \varphi(a_1a_2a_3)
      - k_1(a_1)k_2(a_2,a_3) - k_2(a_1,a_2)k_1(a_3) \\
  & & - k_2(a_1,a_3)k_1(a_2) - k_1(a_1)k_1(a_2)k_1(a_3) \\
& = & \varphi(a_1a_2a_3)
      - \varphi(a_1)\varphi(a_2a_3)
      - \varphi(a_1a_2)\varphi(a_3) \\
  & & - \varphi(a_1a_3)\varphi(a_2)
      + 2 \varphi(a_1)\varphi(a_2)\varphi(a_3).
\end{eqnarray*}
Again we can write this in the Moebius inverted form:
\begin{align*}
k_3(a_1,a_2,a_3)&= \varphi_{\usebox{\NCiiiannew} \,\,\,}[a_1,a_2,a_3]
      - \varphi_{\usebox{\NCiiibnnew} \,\,\,}[a_1,a_2,a_3]\\
&\quad      - \varphi_{\usebox{\NCiiicnnew} \,\,\,}[a_1,a_2,a_3]
      - \varphi_{\usebox{\NCiiidnnew} \,\,\,}[a_1,a_2,a_3]\\
&\quad
      + 2 \varphi_{\usebox{\NCiiiennew} \,\,\,}[a_1,a_2,a_3].
\end{align*}
\end{itemize}

\subsection{Freeness.}
Freeness of subalgebras or random variables is the crucial concept
in free probability theory; it is a non-commutative replacement for
the classical concept of ``independence".
\\
1) Let $\mathcal{ A}_1,\dots,\mathcal{ A}_m \subset \mathcal{ A}$
  be subalgebras with $1 \in \mathcal{ A}_i$ ($i = 1,\dots,m$). The
  subalgebras $\mathcal{ A}_1,\dots,\mathcal{ A}_m$ are called {\it free},
if
  $\varphi(a_1 \cdots a_k) = 0$
for all $k \in {\mathbf{N}}$ and
$a_i \in \mathcal{ A}_{j(i)}$ ($1 \leq j(i) \leq m$)
  such that
$\varphi(a_i) = 0$ for all $i = 1, \dots, k$
    and such that neighbouring elements are from
different subalgebras, i.e., $j(1) \neq j(2) \neq \cdots \neq j(k)$.
\\
2) Let $\mathcal{ X}_1,\dots,\mathcal{ X}_m \subset \mathcal{ A}$
  be subsets of $\mathcal{ A}$. Then
  $\mathcal{ X}_1,\dots,\mathcal{ X}_m$ are called free, if
  $\mathcal{ A}_1,\dots,\mathcal{ A}_m$ are free, where, for $i=1,\dots,m$,
 $\mathcal{ A}_i := \mbox{alg}(1,\mathcal{ X}_i)$
  is the algebra generated by $1$ and $\mathcal{ X}_i$.\\
3) In particular, if the algebras $\mathcal{ A}_i :=
  \mbox{alg}(1,a_i)$ ($i=1,\dots,m$) generated by the elements
  $a_1,\dots,a_m \in \mathcal{
  A}$ are free, then $a_1,\dots,a_m$ are called {\it free} random variables.
If the $*$-algebras generated by the random variables $a_1,\dots,a_m$ are
free, then we call $a_1,\dots,a_m$ {\it $*$-free}.
\\
4) Freeness of random variables can be considered as a rule for
expressing joint moments of free variables in terms of the moments of the
single variables. For example, if $\{a_1,a_2\}$ and $b$ are free,
then the following
identity holds:
\begin{equation} \label{rule}
  \varphi(a_1 b a_2) = \varphi(a_1 a_2) \varphi(b) \: .
\end{equation}
5) The basic fact which shows the relevance of the free cumulants in
connection with freeness is the following characterization of freeness
in terms of cumulants. We will only use this characterization of freeness
in our proofs. Thus, for the purpose of this paper, part (2) of the
following proposition could also be used as the definition of freeness.

\subsection{Proposition \cite{S1}.} Let $(\mathcal{A},\varphi)$ be a
non-commutative probability space and $\mathcal{ A}_1,\dots,\mathcal{
A}_m \subset \mathcal{ A}$ subalgebras. Then the following statements
are equivalent:
\begin{enumerate}
\renewcommand{\labelenumi}{(\arabic{enumi})}
\item The subalgebras $\mathcal{ A}_1,\dots,\mathcal{ A}_m$ are free.
\item For all $n \geq 2$ and all $a_i \in \mathcal{ A}_{j(i)}$
  with $1 \leq j(1), \dots, j(n) \leq m$ we have
$k_n(a_1,\dots,a_n) = 0$
whenever there are some $1 \leq l, k \leq n$ with $j(l) \neq j(k)$.
\end{enumerate}

\section{Main combinatorial result} \label{2}

As mentioned in the Introduction we would like to understand the
behaviour of free cumulants with respect to the multiplicative structure of
our
algebra. The crucial property in a multiplicative
context is associativity. On the level of moments
this just means that we can put brackets arbitrarily; for example we have
$\ff((a_1a_2)a_3)=\ff(a_1(a_2a_3))$.
But the corresponding statement on the level of cumulants is, of course,
not true, i.e. $k_2(a_1a_2,a_3)\not= k_2(a_1,a_2a_3)$ in general.
However, there is still a treatable and nice formula which allows to
deal with free cumulants whose entries are products of random variables.
This formula is the main combinatorial result of this paper
and is presented in
this section.

A special case of that theorem, where only one argument of the cumulant has
the form of a product, appeared in \cite{S2}. However, although our theorem
can be considered as an iteration of that special case, the structure of
that iteration is not clear from the presentation in \cite{S2}. The
main observation here is that this iteration really leads to a beautiful
and useful
combinatorial structure. Our proof will not rely on the special case from
\cite{S2}. It is conceptually much clearer to prove the theorem
directly in its
general form than to do it by iteration.

\subsection{Notation.} The general frame for
our theorem is the following:
Let an increasing sequence of integers be given,
$1 \leq i_1 < i_2 < \cdots <
i_m:=n$
and let $a_1,\dots,a_n$
be random variables.
Then we define new random variables $A_j$ as products of the given $a_i$
according to
$A_j:=a_{i_{j-1}+1} \cdots a_{i_j}$ (where $i_0:=0$).
We want to express a cumulant $k_\tau[A_1,\dots,A_m]$ in terms of
cumulants $k_\pi[a_1,\dots,a_n]$. So let
$\tau$ be a
non-crossing partition of the $m$-tuple $(A_1,\dots,A_m)$.
Then we define
$\hat{\tau} \in NC(a_1,\dots,a_n)$ to be that
partition which we get from $\tau$ by replacing each $A_j$ by
$a_{i_{j-1}+1},\dots,a_{i_j}$, i.e., for $a_i$ being a factor
in $A_k$ and $a_j$ being a factor in $A_l$
we have $a_i\sim_{\hat\tau} a_j$ if and only if $A_k\sim_\tau A_l$.
\\
For example, for
$n=6$ and $A_1:= a_1 a_2$, $A_2 := a_3 a_4 a_5$, $A_3 :=
a_6$ and
\setlength{\unitlength}{0.5cm}
\[ \tau = \{(A_1,A_2),(A_3)\}
   \qquad \entspricht \qquad
  \begin{picture}(3,2)\thicklines
  \put(0,0){\line(0,1){1}}
  \put(0,0){\line(1,0){1}}
  \put(1,0){\line(0,1){1}}
  \put(2,0){\line(0,1){1}}
  \put(0,1.5){\makebox(0,0){$A_1$}}
  \put(1,1.5){\makebox(0,0){$A_2$}}
  \put(2.0,1.5){\makebox(0,0){$A_3$}}
  \end{picture}
\]
we get
\setlength{\unitlength}{0.5cm}
\[ \hat{\tau} = \{(a_1,a_2,a_3,a_4,a_5),(a_6)\}
  \qquad \entspricht \qquad
  \begin{picture}(5,2)\thicklines
  \put(0,0){\line(0,1){1}}
  \put(0,0){\line(1,0){4}}
  \put(1,0){\line(0,1){1}}
  \put(2,0){\line(0,1){1}}
  \put(3,0){\line(0,1){1}}
  \put(4,0){\line(0,1){1}}
  \put(5,0){\line(0,1){1}}
  \put(0,1.5){\makebox(0,0){$a_1$}}
  \put(1,1.5){\makebox(0,0){$a_2$}}
  \put(2,1.5){\makebox(0,0){$a_3$}}
  \put(3,1.5){\makebox(0,0){$a_4$}}
  \put(4,1.5){\makebox(0,0){$a_5$}}
  \put(5,1.5){\makebox(0,0){$a_6$}}
  \end{picture} \: .
\]
Note also in particular, that $\hat\tau=1_n$ if and only if
$\tau=1_m$.

\subsection{Theorem.}
Let $m\in\NN$ and
$1 \leq i_1 < i_2 < \cdots < i_m:=n$ be given.
Consider random variables $a_1,\dots,a_n$
and put
$A_j:=a_{i_{j-1}+1}\cdots a_{i_j}$ for $j=1,\dots,m$ (where $i_0:=0$).
Let $\tau$ be a partition in $NC(A_1,\dots,A_m)$.\\
Then the following equation holds:
\begin{equation} \label{th1_1}
  k_{\tau}[a_1 \cdots a_{i_1}, \dots,
           a_{i_{m-1}+1} \cdots a_{i_m}]
  = \sum_{ \pi \in NC(n) \atop \pi \vee \sigma = \hat{\tau}}
    k_{\pi} [a_1,\dots,a_n] \: ,
\end{equation}
where $\sigma\in NC(n)$ is the partition $\sigma=\{(a_1,\dots,a_{i_1}),
\dots, (a_{i_{m-1}+1},\dots,a_{i_m})\}$.

\

Before we give the proof of our theorem, we want to make clear the structure
of the statement by an example:
\\
\setlength{\unitlength}{0.2cm}
  \newsavebox{\NCiia}
  \savebox{\NCiia}(1,1){
  \thicklines
  \put(0,0){\line(0,1){1}}
  \put(0,0){\line(1,0){1}}
  \put(1,0){\line(0,1){1}}}
\setlength{\unitlength}{0.1cm}
  \newsavebox{\NCiiid}
  \savebox{\NCiiid}(2,2){
  \thicklines
  \put(0,0){\line(0,1){2}}
  \put(0,0){\line(1,0){2}}
  \put(2,0){\line(0,1){2}}
  \put(1,1){\line(0,1){1}}}
\setlength{\unitlength}{0.1cm}
  \newsavebox{\NCiiia}
  \savebox{\NCiiia}(2,2){
  \thicklines
  \put(0,0){\line(0,1){2}}
  \put(0,0){\line(1,0){2}}
  \put(1,0){\line(0,1){2}}
  \put(2,0){\line(0,1){2}}}
\setlength{\unitlength}{0.1cm}
  \newsavebox{\NCiiib}
  \savebox{\NCiiib}(2,2){
  \thicklines
  \put(0,0){\line(0,1){2}}
  \put(1,0){\line(0,1){2}}
  \put(1,0){\line(1,0){1}}
  \put(2,0){\line(0,1){2}}}
\setlength{\unitlength}{0.15cm}
  \newsavebox{\NCiiic}
  \savebox{\NCiiic}(2,2){
  \thicklines
  \put(0,0){\line(0,1){2}}
  \put(0,0){\line(1,0){1}}
  \put(1,0){\line(0,1){2}}
  \put(2,0){\line(0,1){2}}}
\setlength{\unitlength}{0.1cm}
  \newsavebox{\NCiiie}
  \savebox{\NCiiie}(2,2){
  \thicklines
  \put(0,0){\line(0,1){2}}
  \put(1,0){\line(0,1){2}}
  \put(2,0){\line(0,1){2}}}
For $A_1 := a_1 a_2$ and $A_2 := a_3$ we have
$\sigma = \{ (a_1,a_2), (a_3) \} \: \entspricht \usebox{\NCiiic}
\quad$. Consider now
$\tau = 1_2=\{(A_1,A_2)\}$, implying that
$\hat\tau=1_3=\{(a_1,a_2,a_3)\}$.
Then the application of our theorem yields
\begin{eqnarray*}
  k_2(a_1a_2,a_3)
  & = & \sum_{\pi \in NC(3) \atop\pi \vee \sigma = 1_3}
        k_{\pi}[a_1,a_2,a_3] \\
  & = & k_{\usebox{\NCiiia} \,\,\,}[a_1,a_2,a_3]
        + k_{\usebox{\NCiiib} \,\,\,}[a_1,a_2,a_3]
        + k_{\usebox{\NCiiid} \,\,\,}[a_1,a_2,a_3]  \\
  & = & k_3(a_1,a_2,a_3)
        + k_1(a_1) k_2(a_2,a_3) + k_2(a_1,a_3) k_1(a_2) \: ,
\end{eqnarray*}
which is easily seen to be indeed equal to
$k_2(a_1a_2,a_3)=\varphi(a_1 a_2 a_3) -
\varphi(a_1 a_2) \varphi(a_3)$.

\begin{proof}
We show the assertion by induction over
the number $m$ of arguments of the cumulant $k_\tau$.
\\
To begin with, let us study the case when $m=1$. Then we have
$\sigma = \{(a_1, \dots, a_n) \} = 1_n = \hat{\tau}$ and by
the defining relation (1) for the free cumulants our assertion
reduces to
\begin{align*}
k_1(a_1 \cdots  a_n)
  & = \sum_{\pi \in NC(n) \atop \pi \vee 1_n = 1_n}
     k_{\pi} [a_1,\dots,a_n] \\
& = \sum_{\pi \in NC(n)} k_{\pi} [a_1,\dots,a_n] \\
&  = \varphi(a_1 \cdots a_n),
\end{align*}
which is true since $k_1=\ff$.
\\
Let us now make the induction hypothesis that for an integer $m \geq 1$
the theorem is true for all $m' \leq
m$.
\\
We want to show that it also holds for $m+1$. This means that for
$\tau \in NC(m+1)$,
a sequence $1\leq i_1<i_2<\dots<i_{m+1}=:n$, and random variables
$a_1,\dots,a_n$ we have to prove the validity of the following
equation:
\begin{align} \label{th1_proof}
k_{\tau}[A_1,\dots,A_{m+1}]&=
k_{\tau}[a_1 \cdots a_{i_1}, \dots,
           a_{i_m+1} \cdots a_{i_{m+1}}] \nonumber\\
&  = \sum_{\pi \in NC(n) \atop \pi \vee \sigma = \hat{\tau}}
    k_{\pi} [a_1,\dots,a_n],
\end{align}
where $\sigma = \{(a_1,\dots,a_{i_1}), \dots,
(a_{i_m+1},\dots,a_{i_{m+1}})\}$.\\
The proof is divided into two steps. The first one discusses the
case where $\tau \in NC(m+1), \: \tau \neq 1_{m+1}$ and the second
one treats the case where $\tau = 1_{m+1}$.
\begin{description}
\item[\rm Step $1^{\circ}$] The validity of relation
(\ref{th1_proof}) for all $\tau \in NC(m+1)$ except the partition
$1_{m+1}$ is shown as follows: Each such $\tau$ has at least two
blocks, so it can be
written as $\tau = \tau_1 \cup \tau_2$ with $\tau_1$ being a
non-crossing partition of an $s$-tuple $(B_1,\dots,B_s)$ and
$\tau_2$ being a non-crossing partition of a $t$-tuple
$(C_1,\dots,C_t)$ where $(B_1,\dots,B_s)
\cup (C_1,\dots,C_t)=(A_1,\dots,A_{m+1})$ and $s+t = m+1$. With
these definitions, we have
\[ k_{\tau}[A_1,\dots,A_{m+1}]
   = k_{\tau_1} [B_1,\dots,B_s] \: k_{\tau_2} [C_1,\dots,C_t]
                                                             \: .
\]
We will apply now the induction hypothesis on $k_{\tau_1}
[B_1,\dots,B_s]$ and on
$k_{\tau_2} [C_1,\dots,C_t]$.
According to the definition of $A_j$, both $B_k \: (k = 1,\dots,
s)$ and $C_l \: (l = 1,\dots,t)$ are products with factors from
$(a_1,\dots,a_n)$. Put $(b_1,\dots,b_p)$ the tuple containing
all factors of $(B_1,\dots,B_s)$ and $(c_1,\dots,c_q)$ the
tuple consisting of all factors of $(C_1,\dots,C_t)$; this means
$(b_1,\dots,b_p) \cup (c_1,\dots,c_q)=(a_1,\dots,a_n)$ (and
$p+q = n$). We put $\sigma_1 := \sigma \vert_{(b_1,\dots,b_p)}$ and
$\sigma_2 := \sigma\vert_{(c_1,\dots,c_q)}$,
i.e., we have $\sigma=\sigma_1\cup\sigma_2$. Note that $\hat\tau$
factorizes in the same way as $\hat \tau=\hat\tau_1\cup\hat\tau_2$. Then
we get with the
help of our induction hypothesis:
\begin{align*}
k_{\tau}[A_1&,\dots,A_{m+1}]
 =  k_{\tau_1} [B_1,\dots,B_s] \cdot
      k_{\tau_2} [C_1,\dots,C_t] \\
& = \sum_{\pi_1 \in NC(p) \atop
           \pi_1 \vee \sigma_1 = \hat{\tau}_1}
      k_{\pi_1} [b_1,\dots,b_p]
      \cdot
      \sum_{\pi_2 \in NC(q) \atop \pi_2 \vee \sigma_2 = \hat{\tau}_2}
      k_{\pi_2} [c_1,\dots,c_q] \\
& =  \sum_{\pi_1 \in NC(p) \atop \pi_1 \vee \sigma_1 = \hat{\tau}_1}
      \sum_{\pi_2 \in NC(q) \atop \pi_2 \vee \sigma_2 = \hat{\tau}_2}
    k_{\pi_1 \cup \pi_2} [a_1,\dots,a_n] \\
& =  \sum_{\pi \in NC(n) \atop \pi \vee \sigma = \hat{\tau}}
      k_{\pi} [a_1,\dots,a_n] \: .
\end{align*}
\item[\rm Step $2^{\circ}$]
It remains to prove that the equation (\ref{th1_proof}) is also
valid for $\tau = 1_{m+1}$. With (\ref{1.2,3}), we obtain
\begin{align} \label{th1_proof2}
k_{1_{m+1}}[A_1&,\dots,A_{m+1}]
 =  k_{m+1}(A_1,\dots,A_{m+1}) \nonumber \\
& =  \varphi(A_1 \cdots A_{m+1})
      - \sum_{\tau \in NC(m+1) \atop \tau \neq 1_{m+1}}
      k_{\tau}[A_1,\dots,A_{m+1}] \: .
\end{align}
First we transform the sum in (\ref{th1_proof2}) with the
result of step $1^{\circ}$:
\begin{eqnarray*}
\sum_{\tau \in NC(m+1) \atop \tau \neq 1_{m+1}}
k_{\tau}[A_1,\dots,A_{m+1}]
& = & \sum_{\tau \in NC(m+1) \atop \tau \neq 1_{m+1}}
      \sum_{\pi \in NC(n) \atop \pi \vee \sigma = \hat{\tau}}
      k_{\pi} [a_1,\dots,a_n] \\
& = & \sum_{\pi \in NC(n) \atop \pi \vee \sigma \neq 1_n}
      k_{\pi} [a_1,\dots,a_n] \: ,
\end{eqnarray*}
where we used the fact that $\tau=1_{m+1}$ is equivalent
to $\hat\tau=1_n$.
\\
The moment in (\ref{th1_proof2}) can be written as
\[ \varphi(A_1 \cdots A_{m+1})
   = \varphi(a_1 \cdots a_n)
   = \sum_{\pi \in NC(n)} k_{\pi} [a_1,\dots,a_n] \: .
\]
Altogether, we get:
\begin{align*}
k_{m+1}[A_1,\dots,A_{m+1}]
& =  \sum_{\pi \in NC(n)} k_{\pi} [a_1,\dots,a_n]
      - \sum_{\pi \in NC(n) \atop \pi \vee \sigma \neq 1_n}
        k_{\pi} [a_1,\dots,a_n] \\
& =   \sum_{\pi \in NC(n) \atop \pi \vee \sigma = 1_n}
      k_{\pi} [a_1,\dots,a_n] \: .
\end{align*}
\end{description}
\end{proof}

\subsection{Remark.}
In all our applications we will only use the special
case of Theorem 2.2
where $\tau = 1_m$. Then the statement of the theorem is the
following: Consider $m\in\NN$, an increasing sequence
$1 \leq i_1 < i_2 < \cdots < i_m:=n$ and
random variables $a_1,\dots,a_n$.
Put $\sigma:=\{(a_1,\dots,a_{i_1}), \dots,
(a_{i_{m-1}+1},\dots,a_{i_m})\}$. Then we have:
\begin{equation} \label{th1_2}
  k_m[a_1 \cdots a_{i_1}, \dots, a_{i_{m-1}+1} \cdots a_{i_m}]
  = \sum_{\pi \in NC(n) \atop \pi \vee \sigma = 1_n}
    k_{\pi} [a_1,\dots,a_n] \: .
\end{equation}

\

The next proposition, which is from \cite{NS2} (Theorem 1.\/4.\/),
is the basic fact on the multiplication of free random variables. We
want to indicate that our Theorem 2.2 can be used to give
a straightforward and
conceptually simple proof of that statement.

\subsection{Proposition \cite{NS2}} For a positive integer $n$, let
$a_1,\dots,a_n,b_1,\dots,b_n$ be
random variables such that $\{a_1,\dots,a_n\}$ and
$\{b_1,\dots,b_n\}$ are free. Then the following equation holds:
\begin{equation}
k_n(a_1 b_1, \dots, a_n b_n)
   = \sum_{\pi \in NC(n)} k_{\pi} [a_1,\dots,a_n]
                          \: k_{K(\pi)} [b_1,\dots,b_n] \: .
\end{equation}

\begin{proof}
We only give a sketch of the proof.
\\
Applying Theorem 2.2 in
the form mentioned above in Eq. (\ref{th1_2}), we get
\[ k_n(a_1 b_1, \dots, a_n b_n)
   = \sum_{\pi} k_{\pi}[a_1,b_1,\dots,a_n,b_n]
\]
where we have to sum over
\[ \pi \in NC(2n)
   \quad \mbox{with} \quad
   \pi \vee \{(a_1,b_1),\dots, (a_n,b_n)\} = 1_{2n} \: .
\]
Because of the assumption ``$\{a_1,\dots,a_n\},
\{b_1,\dots,b_n\}$ free'' we obtain with Prop.~1.5
that all cumulants vanish with the exception of those
which have only elements from $\{a_1,\dots,a_n\}$ or only elements
from $\{b_1,\dots,b_n\}$
as arguments. This means that all partitions $\pi$ contributing to
the sum must have the form
$\pi = \pi_a \cup \pi_b$ with $\pi_a$ being in
$NC(a_1,\dots,a_n)$ and $\pi_b$ being in $NC(b_1,\dots,b_n)$.
Obviously, for each such $\pi$ we have
\[ k_{\pi}[a_1,b_1,\dots,a_n,b_n]
   = k_{\pi_a}[a_1,\dots,a_n] \: k_{\pi_b}[b_1,\dots,b_n] \: .
\]
One can now convince oneself, that for each $\pi_a\in NC(a_1,\dots,a_n)$
there exists exactly one $\pi_b\in NC(b_1,\dots,b_n)$ such that
$\pi=\pi_a\cup\pi_b$ fulfills the condition
$\pi\vee\{(a_1,b_1),\dots,(a_n,b_n)\}=1_{2n}$ and that this
$\pi_b$ is nothing but the complement of $\pi_a$, i.e., we have to
sum exactly over all $\pi=\pi_a\cup K(\pi_a)$ with $\pi_a\in NC(n)$.
This is the assertion.
\end{proof}

\subsection{Remark} In order to get an idea of the complications
arising in the transition from the tracial to the general non-tracial
case let us consider the following variant of the foregoing
proposition. Let $\{a_1,\dots,a_n\}$ be free from $\{b,c\}$ and consider
the cumulant
$k_n(ba_1c,ba_2c,\dots,ba_nc)$. In the tracial case
this is the same as $k_n(a_1cb,a_2cb,\dots,a_ncb)$ and since
$\{a_1,\dots,a_n\}$ is free from $cb$ our above proposition yields
$$k_n(ba_1c,ba_2c,\dots,ba_nc)=\sum_{\pi\in NC(n)} k_\pi[a_1,\dots,a_n]
k_{K(\pi)}[cb,cb,\dots,cb].$$
In the general situation the structure of the result -- a summation
over $\pi\in NC(n)$ and terms given by a product of cumulants
corresponding to blocks of $\pi$ and
blocks of $K(\pi)$ -- is the same, but
now not always $cb$ appears as argument in the cumulants.
Namely, a careful adaption of our above proof for Prop. 2.4
reveals that we have the following result.

\subsection{Proposition}
For a positive integer $n$
consider random variables $a_1,\dots,a_n,b,c$ such that $\{a_1,\dots,a_n\}$
and $\{b,c\}$ are free. Then we have
\begin{multline}
k_n(ba_1c,ba_2c,\dots,ba_nc)\\=
\sum_{\pi\in NC(n)} k_\pi[a_1,\dots,a_n]
k_{\vert V_r\vert}(bc,bc,\dots,bc)
\prod_{i=1}^{r-1} k_{\vert V_i\vert}(\dashv c,bc,\dots,bc,b\vdash),
\end{multline}
where, for $\pi\in NC(n)$, we have written $K(\pi)=\{V_1,\dots,V_r\}$
such that $V_r$ is the block of $K(\pi)$ containing the last element $n$.
Thus the cumulant corresponding to the block of $K(\pi)$ containing $n$ has
only $bc$ as entries, whereas all the other factors for $K(\pi)$ are of the
form $k_m(\dashv c,bc,\dots,bc,b\vdash)$, which is defined as follows:
$$k_m(\dashv c,b_1,\dots,b_{m-1},b\vdash):=\sum_{\pi\in NC(m+1)\atop
\pi\vee \{(1,m+1),(2),(3),\dots,(m)\}=1_{m+1}}
k_\pi[c,b_1,\dots,b_{m-1},b]$$
for arbitrary random variables $c,b,b_1,\dots,b_{m-1}$.

\subsection{Remarks} 1) Note that the
cumulant $k_m(\dashv c,b_1,\dots,b_{m-1},b\vdash)$ is a cumulant
of order $m$; $c$ and $b$ are to be thought of as the factors of one
argument.
However, in the evaluation of the cumulant one has to take care of the
positions
of $c$ and $b$. For example,
$$k_2(\dashv c,b_1,b\vdash)=\ff(cb_1b)-\ff(cb)\ff(b_1).$$
2) Prop. 2.6 suggests that one might consider also cumulants of the form
\begin{equation}
k^\sigma(a_1,\dots,a_n):=\sum_{\pi\in NC(n)\atop \pi\geq\sigma}
\mu(\pi,1_n)
\ff_\pi[a_1,\dots,a_n]
\end{equation}
for arbitrary $\sigma\in NC(n)$. Note that $k^\sigma$ is not a product
of cumulants like $k_\pi$, but a cumulant of order $\vert \sigma\vert $,
where
each block of $\sigma$ corresponds to an argument given by multiplication
of the corresponding variables $a_i$, but with respectation of the nested
structure of the blocks. If $\sigma$ is of the special form
$\sigma=\{(1,\dots,i_1),\dots,(i_{m-1}+1,\dots,i_m)\}$,
as in Theorem 2.2,
then $k^\sigma$ is nothing but
$$k^{\{(1,\dots,i_1),\dots,(i_{m-1}+1,\dots,i_m)\}}
(a_1,\dots,a_n)=k_m(a_1\cdots a_{i_1},\dots,a_{i_{m-1}+1}\cdots
a_{i_m}),$$ whereas $k_m(\dashv c,b_1,\dots,b_{m-1},b\vdash)$ from Prop. 2.6
reads now as
$$k_m(\dashv c,b_1,\dots,b_{m-1},b\vdash)=
k^{\{(1,m+1),(2),(3),\dots,(m)\}}
(c,b_1,\dots,b_{m-1},b).$$
One should, however, note that the structure of the formula for $k^\sigma$
in terms of moments does not only depend on $\vert \sigma\vert$, but
on the concrete form of $\sigma$ itself. For example, for $\sigma=
\{(1,3),(2),(4)\}$ we have
$$k^\sigma(a_1,b,a_2,c)=\ff(a_1ba_2c)-\ff(a_1ba_2)\ff(c)
-\ff(a_1a_2c)\ff(b)+\ff(a_1a_2)\ff(b)\ff(c),$$
which should be compared with
$$k_3(a,b,c)=\ff(abc)-\ff(ab)\ff(c)-\ff(ac)\ff(b)-\ff(a)\ff(bc)
+2\ff(a)\ff(b)\ff(c).$$
\\
One can generalize Theorem 2.2 for $k^\sigma$ as follows: For $\sigma\in
NC(n)$
and random variables $a_1,\dots,a_n$ we have
\begin{equation}
k^\sigma(a_1,\dots,a_n)=\sum_{\pi\in NC(n)\atop
\pi\vee\sigma=1_n}k_\pi[a_1,\dots,a_n].
\end{equation}
The proof of this statement goes along the same lines as our proof of
Theorem 2.2. We will leave the details to the reader.

\section{Applications to $R$-diagonal elements}  \label{3}

\subsection{Notation {\it (alternating)}.}
Let $a$ be a random variable. A
cumulant $k_{2r}(a_1,\dots,a_{2r})$ with arguments from
$\{a,a^{\ast}\}$ is said to have {\it alternating arguments}, if
there does not exist any $a_i \: (1 \leq i \leq 2r-1)$ with
$a_{i+1} = a_i$. We will also say that the cumulant
$k_{2r}(a_1,\dots,a_{2r})$ is {\it alternating}. Cumulants with
an odd number of arguments will always be considered as not alternating.
\\
Example: The cumulant
$k_6(a,a^{\ast},a,a^{\ast},a,a^{\ast})$ is alternating, whereas \linebreak
$k_8(a,a^{\ast},a^{\ast},a,a,a^{\ast},a,a^{\ast})$
or $k_5(a,a^*,a,a^*,a)$ are not alternating.

\subsection{Definition {\it ($R$-diagonal)}.}
A random variable $a$ is
called {\it $R$-diagonal} if for all $r\in\NN$ we have that
$k_r(a_1,\dots,a_r) = 0$
whenever the arguments $a_1,\dots,a_r\in\{a,a^*\}$
are not alternating in $a$ and $a^*$.

\subsection{Definition {\it (Haar unitary)}.} We call an
element $u$ in a probability space $(\mathcal{A},\varphi)$ {\it
Haar unitary} if it has the following properties:
\begin{enumerate}
\renewcommand{\labelenumi}{(\arabic{enumi}$\mbox{\hspace{0cm}}^{\circ}$)}
\item $u$ is unitary, i.\/e., $u u^{\ast} = 1 = u^{\ast}u$.
\item $\varphi(u^k) = 0 = \varphi(u^{\ast k})$ for $k = 1,2,3,
       \dots $.
\end{enumerate}

\subsection{Remarks.}
1) Due to the relation (\ref{1.2,1}) between moments and free
cumulants, two tuples $(a_1,\dots,a_n)$ and $(b_1,\dots,b_n)$
of random variables have the same joint
distribution if and only if all their cumulants are identical, i.e., if
$k_m(a_{i(1)},\dots,a_{i(m)}) = k_m(b_{i(1)},\dots,b_{i(m)})$ for all
$m\in\NN$ and
all $1 \leq i(1),\dots,i(m) \leq n$.
This implies, of course, that the property ``$R$-diagonality" depends
only on the $*$-distribution of $a$.
\\
2) It was proved in \cite{S3} that a Haar unitary
is $R$-diagonal. Indeed, the examples of the Haar unitary and the
circular element -- which present the two most important non-selfadjoint
distributions in free probability
theory -- provided the motivation for introducing
the class of $R$-diagonal elements as a kind of interpolation between these
two elements.
\\
3) It is clear that all information on the $*$-distribution of an
$R$-diagonal element $a$ is contained in the two sequences of its
alternating cumulants
$\alpha_n:=k_{2n}(a,a^{\ast},a,a^{\ast},\dots,
a,a^{\ast})$ and
$\beta_n:=k_{2n}(a^{\ast},a,a^{\ast},a,\dots,a^{\ast},a)$.
Another useful description of the $*$-distribution of $a$ is given
by the distributions of $aa^*$ and $a^*a$. The next proposition
connects these two descriptions of the $*$-distribution of $a$.
The tracial case -- in which $\alpha_n=\beta_n$ for all $n$ --
was treated in \cite{NS1}, whereas the result in
the general case proves a conjecture, Eq. (5.7), from \cite{NSS}.

\subsection{Proposition.}
Let $a$ be an $R$-diagonal random
variable in a non-commutative probability space $(\mathcal{
A},\varphi)$. Let
\begin{align*}
\alpha_n:&=k_{2n}(a,a^{\ast},a,a^{\ast},\dots,
a,a^{\ast}),\\
\beta_n:&=k_{2n}(a^{\ast},a,a^{\ast},a,\dots,a^{\ast},a)
\end{align*}
be the non-vanishing cumulants of $a$.
Then we have:
\begin{equation}
k_n(aa^*,\dots,aa^*) = \sum_{\pi \in NC(n) \atop \pi = \{V_1,\dots,V_r\}}
              \alpha_{\mid V_1 \mid}
              \beta_{\mid V_2 \mid} \cdots \beta_{\mid V_r \mid},
\end{equation}
where $V_1$ denotes that block of $\pi\in NC(n)$ which contains the
first element 1.

\begin{proof}

Applying Theorem 2.2 in the particular form of Eq. (9) yields
\begin{equation} \label{3_4_proof1}
k_n(aa^{\ast},\dots,aa^{\ast})
  = \sum_{\pi \in NC(2n) \atop \pi \vee \sigma = 1_{2n}}
    k_{\pi} [a,a^{\ast},\dots,a,a^{\ast}]
\end{equation}
with
$$\sigma=\{(a,a^{\ast}),\dots,(a,a^{\ast})\}
\qquad\entspricht\qquad
\{(\mbox{\boldmath$1$},2), \dots,
                  (\mbox{\boldmath$2n-1$},2n)\}.$$

We claim now the following: The partitions $\pi$ which fulfill the
condition $\pi\vee\sigma=1_{2n}$ are exactly those which have the
following properties: the block of $\pi$ which contains the element
$\boldsymbol{1}$
contains also the element $2n$, and, for each $k=1,\dots,n-1$,
the block of $\pi$ which contains the element $2k$
contains also the element \boldmath $2k+1$ \unboldmath.

Since the set of those $\pi\in NC(2n)$ fulfilling the claimed condition
is in canonical bijection with $NC(n)$ and since
$k_\pi[a,a^*,\dots,a,a^*]$ goes
under this bijection to the product
appearing in Eq. (14), this gives directly the assertion.

So it remains to prove the claim.
It is clear that a partition which has the claimed property does
also fulfill $\pi\vee\sigma=1_{2n}$. So we only have to prove the
other direction.

Let $V$ be the block of $\pi$ which contains the element
\boldmath $1$\unboldmath.  Since $a$ is $R$-diagonal the last element
of this block has to be an $a^*$, i.e., an even number, let's say
$2k$. If this would not be $2n$ then this block $V$ would
in $\pi\vee\sigma$ not be
connected to the block containing \boldmath $2k+1$\unboldmath, thus
$\pi\vee\sigma$ would not give $1_{2n}$. Hence $\pi\vee\sigma=1_{2n}$
implies that the block containing the first
element $\boldsymbol{1}$ contains also
the last element $2n$.

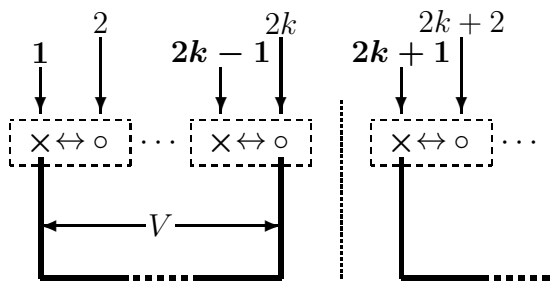
\setlength{\unitlength}{0.4cm}
\[ \begin{picture}(16,9)
   \linethickness{0.6mm}
   \put(1,0){\line(0,1){4}}
   \put(1,0){\line(1,0){2.9}}
   \put(4.1,0){\line(1,0){0.2}}%
   \put(4.5,0){\line(1,0){0.2}}%
   \put(4.9,0){\line(1,0){0.2}}%
   \put(5.3,0){\line(1,0){0.2}}%
   \put(5.7,0){\line(1,0){0.2}}%
   \put(6.1,0){\line(1,0){2.9}}
   \put(9,0){\line(0,1){4}}
   \thicklines
   \put(5,1.75){\makebox(0,0){$V$}}
   \put(4.5,1.75){\vector(-1,0){3.5}}
   \put(5.5,1.75){\vector(1,0){3.5}}
   \linethickness{0.6mm}                                   %
   \put(1,4.5){\makebox(0,0){\boldmath $\times$}}          %
   \put(2,4.5){\makebox(0,0){$\leftrightarrow$}}           %
   \put(3,4.5){\makebox(0,0){$\circ$}}                     %
   \thicklines                                             %
   \put(1,7){\vector(0,-1){1.5}}                           %
   \put(1,7.5){\makebox(0,0){\boldmath$1$\unboldmath}}     %
   \put(3,8){\vector(0,-1){2.5}}                           %
   \put(3,8.5){\makebox(0,0){$2$}}     %
   \thinlines
   \put(0,3.75){\dashbox{0.2}(4,1.5){}}
   \linethickness{0.6mm}
   \put(5,4.5){\makebox(0,0){$\cdots$}}
   \put(7,4.5){\makebox(0,0){\boldmath $\times$}}
   \put(8,4.5){\makebox(0,0){$\leftrightarrow$}}
   \put(9,4.5){\makebox(0,0){$\circ$}}
   \thicklines
   \put(7,7){\vector(0,-1){1.5}}
   \put(7,7.5){\makebox(0,0){\boldmath $2k-1$\unboldmath}}
   \put(9,8){\vector(0,-1){2.5}}
   \put(9,8.5){\makebox(0,0){$2k$}}
   \thinlines
   \put(6,3.75){\dashbox{0.2}(4,1.5){}}
   \thicklines
   \multiput(11,0)(0,0.2){30}{\line(0,1){0.1}}
   \linethickness{0.6mm}
   \put(13,0){\line(0,1){4}}
   \put(13,0){\line(1,0){2.9}}
   \put(16.1,0){\line(1,0){0.2}}%
   \put(16.5,0){\line(1,0){0.2}}%
   \put(16.9,0){\line(1,0){0.2}}%
   \put(17.3,0){\line(1,0){0.2}}%
   \put(17.7,0){\line(1,0){0.2}}%
   \put(13,4.5){\makebox(0,0){\boldmath $\times$}}
   \put(14,4.5){\makebox(0,0){$\leftrightarrow$}}
   \put(15,4.5){\makebox(0,0){$\circ$}}
   \thicklines
   \put(13,7){\vector(0,-1){1.5}}
   \put(13,7.5){\makebox(0,0){\boldmath$2k+1$}}
   \put(15,8){\vector(0,-1){2.5}}
   \put(15,8.5){\makebox(0,0){$2k+2$}}
   \thinlines
   \put(12,3.75){\dashbox{0.2}(4,1.5){}}
   \thicklines
   \put(17,4.5){\makebox(0,0){$\cdots$}}
   \thinlines
   \end{picture}
\]

Now fix a $k=1,\dots,n-1$ and let $V$ be the block of $\pi$ containing
the element $2k$. Assume that $V$ does not contain the element
\boldmath $2k+1$\unboldmath. Then there are two possibilities: Either
$2k$ is not the last element in $V$, i.e. there exists a next
element in $V$, which is necessarily of the form $\boldsymbol{2l+1}$ with
$l>k$
...

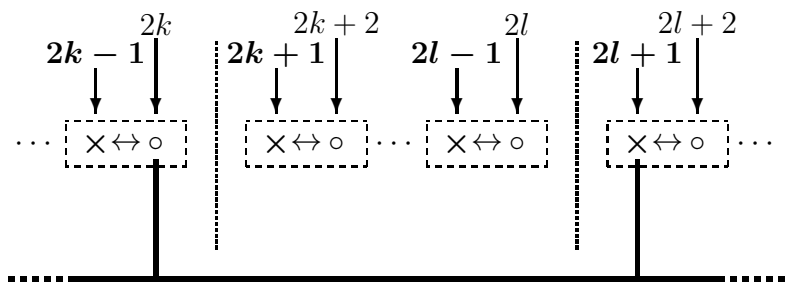
\setlength{\unitlength}{0.4cm}
\[ \begin{picture}(34,9)
   \linethickness{0.6mm}
   \put(4.1,0){\line(1,0){0.2}}%
   \put(4.5,0){\line(1,0){0.2}}%
   \put(4.9,0){\line(1,0){0.2}}%
   \put(5.3,0){\line(1,0){0.2}}%
   \put(5.7,0){\line(1,0){0.2}}%
   \put(6.1,0){\line(1,0){2.9}}
   \put(9,0){\line(0,1){4}}
   \put(9,0){\line(1,0){16}}
   \thicklines
   \linethickness{0.6mm}                                   %
   \thinlines
   \linethickness{0.6mm}
   \put(5,4.5){\makebox(0,0){$\cdots$}}
   \put(7,4.5){\makebox(0,0){\boldmath $\times$}}
   \put(8,4.5){\makebox(0,0){$\leftrightarrow$}}
   \put(9,4.5){\makebox(0,0){$\circ$}}
   \thicklines
   \put(7,7){\vector(0,-1){1.5}}
   \put(7,7.5){\makebox(0,0){\boldmath $2k-1$\unboldmath}}
   \put(9,8){\vector(0,-1){2.5}}
   \put(9,8.5){\makebox(0,0){$2k$}}
   \thinlines
   \put(6,3.75){\dashbox{0.2}(4,1.5){}}
   \thicklines
   \multiput(11,1)(0,0.2){35}{\line(0,1){0.1}}
   \linethickness{0.6mm}
   \linethickness{0.6mm}
   \put(13,4.5){\makebox(0,0){\boldmath $\times$}}
   \put(14,4.5){\makebox(0,0){$\leftrightarrow$}}
   \put(15,4.5){\makebox(0,0){$\circ$}}
   \thicklines
   \put(13,7){\vector(0,-1){1.5}}
   \put(13,7.5){\makebox(0,0){\boldmath$2k+1$}}
   \put(15,8){\vector(0,-1){2.5}}
   \put(15,8.5){\makebox(0,0){$2k+2$}}
   \thinlines
   \put(12,3.75){\dashbox{0.2}(4,1.5){}}
   \put(17,4.5){\makebox(0,0){$\cdots$}}
   \thicklines
   \linethickness{0.6mm}
   \put(19,4.5){\makebox(0,0){\boldmath $\times$}}
   \put(20,4.5){\makebox(0,0){$\leftrightarrow$}}
   \put(21,4.5){\makebox(0,0){$\circ$}}
   \thicklines
   \put(19,7){\vector(0,-1){1.5}}
   \put(19,7.5){\makebox(0,0){\boldmath $2l-1$}}
   \put(21,8){\vector(0,-1){2.5}}
   \put(21,8.5){\makebox(0,0){$2l$}}
   \thinlines
   \put(18,3.75){\dashbox{0.2}(4,1.5){}}
   \thicklines
   \multiput(23,1)(0,0.2){35}{\line(0,1){0.1}}
   \linethickness{0.6mm}
   \put(25,0){\line(0,1){4}}
   \put(25,0){\line(1,0){2.9}}
   \put(28.1,0){\line(1,0){0.2}}%
   \put(28.5,0){\line(1,0){0.2}}%
   \put(28.9,0){\line(1,0){0.2}}%
   \put(29.3,0){\line(1,0){0.2}}%
   \put(29.7,0){\line(1,0){0.2}}%
   \thicklines
   \linethickness{0.6mm}
   \put(25,4.5){\makebox(0,0){\boldmath $\times$}}
   \put(26,4.5){\makebox(0,0){$\leftrightarrow$}}
   \put(27,4.5){\makebox(0,0){$\circ$}}
   \thicklines
   \put(25,7){\vector(0,-1){1.5}}
   \put(25,7.5){\makebox(0,0){\boldmath $2l+1$}}
   \put(27,8){\vector(0,-1){2.5}}
   \put(27,8.5){\makebox(0,0){$2l+2$}}
   \thinlines
   \put(24,3.75){\dashbox{0.2}(4,1.5){}}
   \linethickness{0.6mm}
   \put(29,4.5){\makebox(0,0){$\cdots$}}
   \end{picture}
\]

... or $2k$ is the last element in $V$. In this case the first element
of $V$ is of the form $\boldsymbol{2l+1}$ with $0\leq l\leq k-1$.

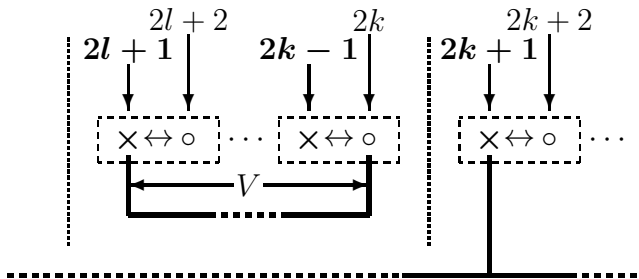
\setlength{\unitlength}{0.4cm}
\[ \begin{picture}(34,9)
   \linethickness{0.6mm}
   \multiput(9,0)(0.4,0){33}{\line(1,0){0.2}}
\thicklines
   \multiput(11,1)(0,0.2){35}{\line(0,1){0.1}}
   \linethickness{0.6mm}
   \put(13,4.5){\makebox(0,0){\boldmath $\times$}}
   \put(14,4.5){\makebox(0,0){$\leftrightarrow$}}
   \put(15,4.5){\makebox(0,0){$\circ$}}
   \thicklines
   \put(13,7){\vector(0,-1){1.5}}
   \put(13,7.5){\makebox(0,0){\boldmath$2l+1$}}
   \put(15,8){\vector(0,-1){2.5}}
   \put(15,8.5){\makebox(0,0){$2l+2$}}
   \thinlines
   \put(12,3.75){\dashbox{0.2}(4,1.5){}}
   \put(17,4.5){\makebox(0,0){$\cdots$}}
   \thicklines
   \linethickness{0.6mm}
   \put(19,4.5){\makebox(0,0){\boldmath $\times$}}
   \put(20,4.5){\makebox(0,0){$\leftrightarrow$}}
   \put(21,4.5){\makebox(0,0){$\circ$}}
   \thicklines
   \put(19,7){\vector(0,-1){1.5}}
   \put(19,7.5){\makebox(0,0){\boldmath $2k-1$}}
   \put(21,8){\vector(0,-1){2.5}}
   \put(21,8.5){\makebox(0,0){$2k$}}
   \thinlines
   \put(18,3.75){\dashbox{0.2}(4,1.5){}}
   \thicklines
   \multiput(23,1)(0,0.2){35}{\line(0,1){0.1}}
   \linethickness{0.6mm}
   \put(25,0){\line(0,1){4}}
   \put(22.1,0){\line(1,0){5.8}}
   \put(28.1,0){\line(1,0){0.2}}%
   \put(28.5,0){\line(1,0){0.2}}%
   \put(28.9,0){\line(1,0){0.2}}%
   \put(29.3,0){\line(1,0){0.2}}%
   \put(29.7,0){\line(1,0){0.2}}%
   \thicklines
   \linethickness{0.6mm}
   \put(25,4.5){\makebox(0,0){\boldmath $\times$}}
   \put(26,4.5){\makebox(0,0){$\leftrightarrow$}}
   \put(27,4.5){\makebox(0,0){$\circ$}}
   \thicklines
   \put(25,7){\vector(0,-1){1.5}}
   \put(25,7.5){\makebox(0,0){\boldmath $2k+1$}}
   \put(27,8){\vector(0,-1){2.5}}
   \put(27,8.5){\makebox(0,0){$2k+2$}}
   \thinlines
   \put(24,3.75){\dashbox{0.2}(4,1.5){}}
   \thicklines
   \put(29,4.5){\makebox(0,0){$\cdots$}}
   \put(17,3){\makebox(0,0){$V$}}
   \put(16.5,3){\vector(-1,0){3.5}}
   \put(17.5,3){\vector(1,0){3.5}}
   \linethickness{0.6mm}
   \put(13,2){\line(0,1){2}}
   \put(13,2){\line(1,0){2.9}}
   \put(16.1,2){\line(1,0){0.2}}%
   \put(16.5,2){\line(1,0){0.2}}%
   \put(16.9,2){\line(1,0){0.2}}%
   \put(17.3,2){\line(1,0){0.2}}%
   \put(17.7,2){\line(1,0){0.2}}%
   \put(18.1,2){\line(1,0){2.9}}
   \put(21,2){\line(0,1){2}}
   \end{picture}
\]
In both cases the block $V$ gets not connected with $\boldsymbol{2k+1}$
in $\pi\vee\sigma$, thus this cannot give $1_{2n}$. Hence the condition
$\pi\vee\sigma=1_{2n}$ forces $2k$ and $\boldsymbol{2k+1}$ to lie in the
same block.
This proves our claim and hence the assertion.
\end{proof}

We are now going to prove a fundamental characterization of
$R$-diagonal elements as those random variables whose
$*$-distribution remains invariant under the multiplication with
a free Haar unitary.
This theorem has been proven in \cite{NS3}
in the case when $\varphi$ is a trace. The treatment there used
some ad hoc combinatorics. In contrast to this, our approach here is
more straightforward and conceptually clearer. Another proof of the
general form of the theorem, relying on Fock space techniques,
will appear in \cite{NSS}.
The main step in the proof of the theorem --
the one in which we will use our combinatorial
Theorem 2.2 --
is the following
proposition. This appeared also, for the tracial case, in \cite{NS1}.

\subsection{Proposition.}
Let $a$ and $x$ be elements in a
probability space $(\mathcal{A},\varphi)$ with $a$ being $R$-diagonal
and such that $\{a,a^{\ast}\}$ and $\{x,x^{\ast}\}$ are free. Then
$ax$ is $R$-diagonal.

\begin{proof}
We examine a cumulant
$k_r(a_1a_2,\dots,a_{2r-1}a_{2r})$ with $a_{2i-1}a_{2i} \in \{
ax,x^{\ast}a^{\ast} \}$ for $i \in \{ 1,\dots,r \}$.\\[0.5ex]
According to the definition of $R$-diagonality we have to show that
this cumulant vanishes in the following two cases:
\begin{enumerate}
\renewcommand{\labelenumi}{(\arabic{enumi}$\mbox{\hspace{0cm}}^{\circ}$)}
\item $r$ is odd.
\item There exists
  at least one $s \: (1 \leq s \leq r-1)$ such that
  $a_{2s-1}a_{2s} = a_{2s+1}a_{2s+2}$.
\end{enumerate}

By Theorem 2.2, we have
\begin{equation}
k_r(a_1a_2,\dots,a_{2r-1}a_{2r})
   = \sum_{\pi \in NC(2r) \atop \pi \vee \sigma = 1_{2r}}
     k_{\pi}[a_1,a_2,\dots,a_{2r-1},a_{2r}] \: ,
\end{equation}
where
$\sigma = \{(a_1,a_2),\dots,(a_{2r-1},a_{2r})\}$.

The fact that $a$ and $x$ are $*$-free
implies, by Prop.~1.5, that only such partitions $\pi
\in NC(2r)$ contribute to the sum each of whose blocks contains elements
only from $\{a,a^{\ast}\}$ or only from $\{x,x^{\ast}\}$.

{\rm Case ($1^{\circ}$):} As there is at least one block
of $\pi$ containing a different number of
elements $a$ and $a^{\ast}$, $k_\pi$ vanishes always.
So there are no partitions $\pi$ contributing to the
sum in (16) which consequently vanishes.

{\rm Case ($2^{\circ}$):} We assume that there exists an
$s \in \{1,\dots,r-1\}$ such that $a_{2s-1}a_{2s} =
a_{2s+1}a_{2s+2}$. Since with $a$ also $a^*$ is $R$-diagonal, it
suffices to consider the case
where $a_{2s-1}a_{2s} = a_{2s+1}a_{2s+2}=ax$, i.\/e.,
$a_{2s-1} = a_{2s+1}=a$ and $a_{2s} = a_{2s+2}=x$.
\\
Let $V$ be the block containing $a_{2s+1}$. We have to examine
two situations:
\begin{enumerate}
\renewcommand{\labelenumi}{\Alph{enumi}.}
\item On the one hand, it might happen that $a_{2s+1}$ is the
first element in the block $V$. This can be sketched in the following way:
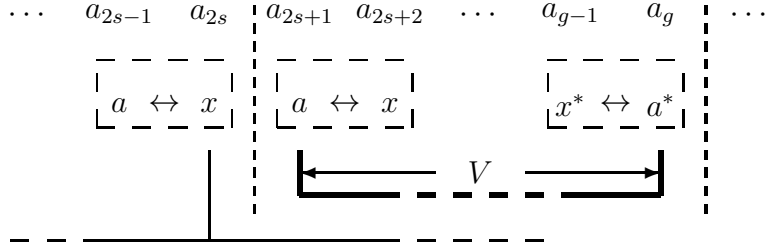
\setlength{\unitlength}{1.2cm}
\[ \begin{picture}(10,3)
   \thicklines
   \put(3,0){\line(0,1){1}}
   \put(1.6,0){\line(1,0){3.5}}
   \put(1.2,0){\line(1,0){0.2}}
   \put(0.8,0){\line(1,0){0.2}}
   \put(5.3,0){\line(1,0){0.2}}
   \put(5.7,0){\line(1,0){0.2}}
   \put(6.1,0){\line(1,0){0.2}}
   \put(6.5,0){\line(1,0){0.2}}
   \linethickness{0.6mm}
   \put(4,0.5){\line(0,1){0.5}}
   \put(4,0.5){\line(1,0){1.1}}
   \put(5.3,0.5){\line(1,0){0.2}}
   \put(5.7,0.5){\line(1,0){0.2}}
   \put(6.1,0.5){\line(1,0){0.2}}
   \put(6.5,0.5){\line(1,0){0.2}}
   \put(6.9,0.5){\line(1,0){1.1}}
   \put(8,0.5){\line(0,1){0.5}}
   \put(2,1.5){\makebox(0,0){$a$}}
   \put(2.5,1.5){\makebox(0,0){$\leftrightarrow$}}
   \put(3,1.5){\makebox(0,0){$x$}}
   \put(4,1.5){\makebox(0,0){$a$}}
   \put(4.5,1.5){\makebox(0,0){$\leftrightarrow$}}
   \put(5,1.5){\makebox(0,0){$x$}}
   \put(7,1.5){\makebox(0,0){$x^{\ast}$}}
   \put(7.5,1.5){\makebox(0,0){$\leftrightarrow$}}
   \put(8,1.5){\makebox(0,0){$a^{\ast}$}}
   \put(1,2.5){\makebox(0,0){$\cdots$}}
   \put(2,2.5){\makebox(0,0){$a_{2s-1}$}}
   \put(3,2.5){\makebox(0,0){$a_{2s}$}}
   \thicklines
   \multiput(3.5,0.3)(0,0.2){12}{\line(0,1){0.1}} %
   \linethickness{0.6mm}
   \put(4,2.5){\makebox(0,0){$a_{2s+1}$}}
   \put(5,2.5){\makebox(0,0){$a_{2s+2}$}}
   \put(6,2.5){\makebox(0,0){$\cdots$}}
   \put(7,2.5){\makebox(0,0){$a_{g-1}$}}
   \put(8,2.5){\makebox(0,0){$a_g$}}
   \thicklines
   \multiput(8.5,0.3)(0,0.2){12}{\line(0,1){0.1}} %
   \linethickness{0.6mm}
   \put(9,2.5){\makebox(0,0){$\cdots$}}
   \thinlines
   \put(1.75,1.25){\dashbox{0.2}(1.5,0.75){}}
   \put(3.75,1.25){\dashbox{0.2}(1.5,0.75){}}
   \put(6.75,1.25){\dashbox{0.2}(1.5,0.75){}}
   \thicklines
   \put(6,0.75){\makebox(0,0){$V$}}
   \put(6.5,0.75){\vector(1,0){1.5}}
   \put(5.5,0.75){\vector(-1,0){1.5}}
   \end{picture}
\]
In this case the block $V$ is not connected with $a_{2s}$ in
$\pi\vee\sigma$, thus the latter cannot be equal to $1_{2n}$.
\item On the other hand, it can happen that $a_{2s+1}$ is not
the first element of $V$. Because $a$ is $R$-diagonal, the preceeding
element must be an $a^*$.
\setlength{\unitlength}{1.2cm}
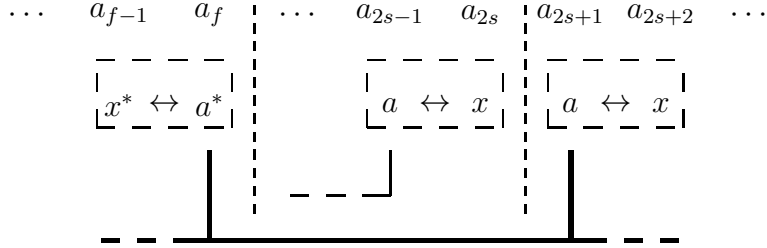
\[ \begin{picture}(10,3)
   \linethickness{0.6mm}
   \put(3,0){\line(0,1){1}}
   \put(2.2,0){\line(1,0){0.2}}
   \put(1.8,0){\line(1,0){0.2}}
   \put(2.6,0){\line(1,0){4.8}}
   \put(7,0){\line(0,1){1}}
   \put(7.6,0){\line(1,0){0.2}}
   \put(8.0,0){\line(1,0){0.2}}
   \thicklines
   \put(3.9,0.5){\line(1,0){0.2}}
   \put(4.3,0.5){\line(1,0){0.2}}
   \put(4.7,0.5){\line(1,0){0.3}}
   \put(5,0.5){\line(0,1){0.5}}
   \linethickness{0.6mm}
   \put(2,1.5){\makebox(0,0){$x^{\ast}$}}
   \put(2.5,1.5){\makebox(0,0){$\leftrightarrow$}}
   \put(3,1.5){\makebox(0,0){$a^{\ast}$}}
   \put(5,1.5){\makebox(0,0){$a$}}
   \put(5.5,1.5){\makebox(0,0){$\leftrightarrow$}}
   \put(6,1.5){\makebox(0,0){$x$}}
   \put(7,1.5){\makebox(0,0){$a$}}
   \put(7.5,1.5){\makebox(0,0){$\leftrightarrow$}}
   \put(8,1.5){\makebox(0,0){$x$}}
   \put(1,2.5){\makebox(0,0){$\cdots$}}
   \put(2,2.5){\makebox(0,0){$a_{f-1}$}}
   \put(3,2.5){\makebox(0,0){$a_f$}}
   \thicklines
   \multiput(3.5,0.3)(0,0.2){12}{\line(0,1){0.1}} %
   \linethickness{0.6mm}
   \put(4,2.5){\makebox(0,0){$\cdots$}}
   \put(5,2.5){\makebox(0,0){$a_{2s-1}$}}
   \put(6,2.5){\makebox(0,0){$a_{2s}$}}
   \thicklines
   \multiput(6.5,0.3)(0,0.2){12}{\line(0,1){0.1}} %
   \linethickness{0.6mm}
   \put(7,2.5){\makebox(0,0){$a_{2s+1}$}}
   \put(8,2.5){\makebox(0,0){$a_{2s+2}$}}
   \put(9,2.5){\makebox(0,0){$\cdots$}}
   \thinlines
   \put(1.75,1.25){\dashbox{0.2}(1.5,0.75){}}
   \put(4.75,1.25){\dashbox{0.2}(1.5,0.75){}}
   \put(6.75,1.25){\dashbox{0.2}(1.5,0.75){}}
   \end{picture}
\]
But then $V$ will again not be connected to $a_{2s}$ in $\pi\vee\sigma$.
Thus again $\pi\vee\sigma$ cannot be equal to $1_{2n}$.
\end{enumerate}
As in both cases we do not find any partition contributing to the
investigated sum in Eq. (16) this has to vanish.
\end{proof}

\subsection{Theorem.} Let $x$ be an element in a
non-commutative probability space $(\mathcal{A},\varphi)$.
Furthermore, let $u$ be a Haar unitary in
$(\mathcal{A},\varphi)$ such that $\{u,u^{\ast}\}$ and $\{x,x^{\ast}\}$
are free.
Then $x$ is $R$-diagonal if and only if $(x,x^{\ast})$ has the
same joint distribution as $(ux,x^{\ast}u^{\ast})$:
\[\mbox{$x$ $R$-diagonal}
   \quad \Longleftrightarrow \quad
   \mu_{x,x^{\ast}} =
   \mu_{ux,x^{\ast}u^{\ast}} \: . \]

\begin{proof}
$\Longrightarrow$:
In order to show that the joint distributions of $(x,x^{\ast})$
and $(ux,x^{\ast}u^{\ast})$ are identical, we have to prove
according to the Remark 3.4(1) that
$k_m(b_1,\dots,b_m) = k_m(c_1,\dots,c_m)$ for all $m\in\NN$,
$b_i \in \{x,x^{\ast}\}$ and
\[ c_i = \left\{ \begin{array}{l@{\quad\mbox{for}\quad b_i=}l}
                          ux & x \\ x^{\ast}u^{\ast} & x^{\ast}
                 \end{array}
         \right. .
\]
In the cases when $m$ is odd or when with even $m$ the elements
$b_1,\dots,b_m$ do not alternate, the cumulant
$k_m(b_1,\dots,b_m)$ vanishes because of the $R$-diagonality
of $x$. By Prop. 3.6 and the fact that $u$ is $R$-diagonal,
we get that $ux$ is $R$-diagonal, too, and
therefore $k_m(c_1,\dots,c_m)$ also vanishes.\\
Hence we have to consider the case where the arguments
$b_1,\dots,b_m$ alternate (which implies alternating arguments
$c_1,\dots,c_m$).\\
We inductively show the validity of
\[ k_{2r}(x,x^{\ast},\dots,x,x^{\ast})
   = k_{2r}(ux,x^{\ast}u^{\ast},\dots,ux,x^{\ast}u^{\ast})
\]
and
\[ k_{2r}(x^{\ast},x,\dots,x^{\ast},x)
   = k_{2r}(x^{\ast}u^{\ast},ux,\dots,x^{\ast}u^{\ast},ux)
\]
for any natural $r$.

First, consider $r=1$. On one hand, the equation
\[ k_2(ux,x^{\ast}u^{\ast})
   = \varphi(uxx^{\ast}u^{\ast}) - k_1(ux) k_1(x^{\ast}u^{\ast})
\]
holds by definition
of $k_2$. With both cumulants $k_1(ux)$ and
$k_1(x^{\ast}u^{\ast})$ vanishing because of the $R$-diagonality
of $ux$ the second term of the sum is equal to zero.\\
Since $\{u,u^{\ast}\}$ and $\{x,x^{\ast}\}$ are assumed to be
free, we can write the moment with the help of formula
(5) as
$$\varphi(uxx^{\ast}u^{\ast})
  = \varphi(uu^{\ast}) \varphi(xx^{\ast})=\ff(xx^*).$$
So we get
$k_2(ux,x^{\ast}u^{\ast}) = \varphi(xx^{\ast})$.
\\
On the other hand, with $x$ being $R$-diagonal we obtain
\[ k_2(x,x^{\ast})
   = \varphi(xx^{\ast}) - k_1(x)k_1(x^{\ast})
   = \varphi(xx^{\ast})
   = k_2(ux,x^{\ast}u^{\ast}) \: .
\]

Induction hypothesis: Assume the following to be
true for any $r' < r \: (r \geq 2)$:
\begin{align*}
k_{2r'}(x,x^{\ast},\dots,x,x^{\ast})
   &= k_{2r'}(ux,x^{\ast}u^{\ast},\dots,ux,x^{\ast}u^{\ast}) \\
k_{2r'}(x^{\ast},x,\dots,x^{\ast},x)
   &= k_{2r'}(x^{\ast}u^{\ast},ux,\dots,x^{\ast}u^{\ast},ux).
\end{align*}
We have to show the validity of these equations for $r'=r$.
It suffices to consider the first equation.
\\
According to definition of the free cumulants we have
\begin{eqnarray*}
  \lefteqn{
  k_{2r}(ux,x^{\ast}u^{\ast},\dots,ux,x^{\ast}u^{\ast}) } \\
  & = & \varphi(ux x^{\ast}u^{\ast} \cdots ux x^{\ast}u^{\ast})
     - \sum_{\pi \in NC(2r) \atop \pi \neq 1_{2r}}
     k_{\pi}[ux,x^{\ast}u^{\ast},\dots,ux,x^{\ast}u^{\ast}] \: .
\end{eqnarray*}
Because of the freeness of $\{u,u^{\ast}\}$ and $\{x,x^{\ast}\}$
and with the help of (5) we get
\[ \varphi(ux x^{\ast}u^{\ast} ux \cdots x^{\ast}u^{\ast} ux
                                           x^{\ast}u^{\ast})
  = \varphi(u [x x^{\ast}]^r u^{\ast})
  = \varphi(u u^{\ast}) \varphi([x x^{\ast}]^r)
  = \varphi([x x^{\ast}]^r) \: .
\]
It follows that
\[ k_{2r}(ux,x^{\ast}u^{\ast},\dots,ux,x^{\ast}u^{\ast})
   = \varphi([x x^{\ast}]^r)
     - \sum_{\pi \in NC(2r) \atop \pi \neq 1_{2r}}
     k_{\pi}[ux,x^{\ast}u^{\ast},\dots,ux,x^{\ast}u^{\ast}] \: .
\]
The only partitions $\pi \in NC(2r) , \: \pi \neq 1_{2r}$
contributing in
the foregoing sum are those where all blocks are alternating in
$ux$ and $x^*u^*$.
According to our induction hypothesis, we can then replace in all
blocks the element $ux$ by $x$
and the element $x^{\ast}u^{\ast}$ by $x^{\ast}$. So we finally
obtain
\[ k_{2r}(ux,x^{\ast}u^{\ast},\dots,ux,x^{\ast}u^{\ast})
   = k_{2r}(x,x^{\ast},\dots,x,x^{\ast}) \: .
\]

$\Longleftarrow$:
We assume that $\mu_{x,x^{\ast}} = \mu_{ux,x^{\ast}u^{\ast}}$.
As, by Prop. 3.6, $ux$ is $R$-diagonal,
$x$ is $R$-diagonal, too.
\end{proof}

\subsection{Remark}
Prop. 3.6 implies in particular that the product of two free
$R$-diagonal elements is $R$-diagonal again. This raises the
question how the alternating cumulants of the product are given
in terms of the alternating cumulants of the factors. This is
answered in the next proposition. In the tracial case this
reproduces a result of \cite{NS1}, whereas in the general case
this proves the conjecture (5.8) from \cite{NSS}.

\subsection{Proposition.}
Let $a$ and $b$ be $R$-diagonal
random variables such that $\{a,a^{\ast}\}$ is free from
$\{b,b^{\ast}\}$. Furthermore, put
\begin{eqnarray*}
\alpha_n &:=& k_{2n}(a,a^{\ast},a,a^{\ast},\dots,a,a^{\ast})
                                                           \: ,\\
\beta_n  &:=& k_{2n}(a^{\ast},a,a^{\ast},a,\dots,a^{\ast},a)
                                                           \: ,\\
\gamma_n &:=& k_{2n}(b,b^{\ast},b,b^{\ast},\dots,b,b^{\ast}).
\end{eqnarray*}
Then $ab$ is $R$-diagonal and the alternating cumulants of
$ab$ are given by
\begin{multline} \label{3_5}
k_{2n}(ab,b^{\ast}a^{\ast},\dots,ab,b^{\ast}a^{\ast})
\\  = \sum_{{\pi=\pi_a\cup\pi_b \in NC(2n) \atop
 \pi_a = \{V_1,\dots,V_k\}\in NC(1,3,\dots,2n-1)}\atop
\pi_b=\{V_1',\dots,V_l'\}\in NC(2,4,\dots,2n)}
          \alpha_{\mid V_1 \mid}
          \beta_{\mid V_2 \mid} \cdots \beta_{\mid V_k \mid}
          \gamma_{\mid V_1'\mid} \cdots \gamma_{\mid V_l' \mid}
                                                             \: ,
\end{multline}
where $V_1$ is that block of $\pi$ which contains the first element 1.

\begin{proof}
$R$-diagonality of $ab$ is clear by Prop. 3.6. So we only have to
prove Eq.~(17).

By Theorem 2.2, we get
\begin{equation} \label{3_5_proof1}
k_{2n}(ab,b^{\ast}a^{\ast},\ldots,ab,b^{\ast}a^{\ast})
= \sum_{\pi\in NC(4n)\atop \pi\vee\sigma=1_{4n}}
  k_{\pi} [a,b,b^{\ast},a^{\ast},\ldots,a,b,b^{\ast},a^{\ast}]
                                                             \: ,
\end{equation}
where $\sigma = \{(a,b),(b^{\ast},a^{\ast}),\ldots,
                  (a,b),(b^{\ast},a^{\ast})\}$.
Since $\{a,a^{\ast}\}$ and $\{b,b^{\ast}\}$ are assumed to be
free, we also know, by Prop.~1.5, that for a contributing partition
$\pi$ each block has to contain components only
from $\{a,a^{\ast}\}$ or only from $\{b,b^{\ast}\}$.\\
As in the proof of Prop.~3.5 one can show that
the requirement $\pi\vee\sigma=1_{4n}$ is equivalent to the following
properties of $\pi$: The block containing 1 must also contain $4n$ and,
for each $k=1,...,2n-1$, the block containing $2k$ must also contain $2k+1$.
(This couples always $b$ with $b^*$ and $a^*$ with $a$, so it is
compatible with the $*$-freeness between $a$ and $b$.)
The set of partitions in $NC(4n)$
fulfilling these properties is in canonical bijection
with $NC(2n)$. Furthermore we have to take care of the fact that each
block of $\pi\in NC(4n)$ contains either only elements from $\{a,a^*\}$
or only elements from $\{b,b^*\}$. For the image of $\pi$ in $NC(2n)$
this means that it splits into blocks living on the odd numbers and
blocks living on the even numbers. Furthermore, under these identifications
the quantity
$k_{\pi} [a,b,b^{\ast},a^{\ast},\ldots,a,b,b^{\ast},a^{\ast}]$ goes
over to the expression as appearing in our assertion (17).
\end{proof}

\subsection{Remark.} According to Prop.~3.6 multiplication preserves
$R$-diago\-na\-lity if the factors are free. Haagerup and Larsen \cite{HL,L}
showed that,
in the tracial case, the same statement is also true for the
other extreme relation between the factors, namely if they are the same --
i.e., powers of $R$-diagonal elements are also $R$-diagonal. The proof
of Haagerup and Larsen relied on special realizations of $R$-diagonal
elements. Here we will give a short combinatorial proof of that statement.
In particular, our proof will -- in comparison with the proof of
Prop. 3.6 -- also illuminate the relation between
the statements ``$a_1,\dots,a_r$ $R$-diagonal and free implies
$a_1\cdots a_r$ $R$-diagonal" and ``$a$ $R$-diagonal implies $a^r$
$R$-diagonal". Furthermore, in contrast to the approach of \cite{HL,L},
our proof extends without problems to the non-tracial situation.

\subsection{Proposition.} Let $a$ be an $R$-diagonal element
and let $r$ be a positive integer. Then $a^r$ is $R$-diagonal, too.

\begin{proof}
For notational convenience we deal with the case $r=3$.
General $r$ can be treated
analogously.\\
The cumulants which we must have a look at are
$k_n(b_1,\dots,b_n)$ with arguments $b_i$ from $\{a^3,
(a^3)^{\ast}\} \: (i = 1,\dots,n)$. We write $b_i = b_{i,1}
b_{i,2} b_{i,3}$ with $b_{i,1} = b_{i,2} = b_{i,3} \in
\{a,a^{\ast}\}$. According to the definition of
$R$-diagonality we have to show that for any $n \geq 1$ the
cumulant \linebreak
$k_n(b_{1,1} b_{1,2} b_{1,3}, \dots, b_{n,1} b_{n,2}
b_{n,3})$ vanishes if (at least) one of the following things
happens:
\begin{enumerate}
\renewcommand{\labelenumi}{(\arabic{enumi}$\mbox{\hspace{0cm}}^{\circ}$)}
\item There exists an $s \in \{1,
      \dots,n-1\}$ with $b_s=b_{s+1}$.
\item $n$ is odd.
\end{enumerate}

Theorem 2.2 yields
\begin{multline*}
 k_n(b_{1,1} b_{1,2} b_{1,3}, \dots, b_{n,1} b_{n,2} b_{n,3})
  \\  = \sum_{\pi \in NC(3n) \atop \pi \vee \sigma = 1_{3n}}
         k_{\pi} [b_{1,1},b_{1,2},b_{1,3}, \dots,
                  b_{n,1},b_{n,2},b_{n,3}] \: ,
\end{multline*}
where
$\sigma:=\{(b_{1,1},b_{1,2},b_{1,3}), \dots,
(b_{n,1},b_{n,2},b_{n,3})\}$.
The $R$-diagonality of $a$ implies that a partition
$\pi$ gives a non-vanishing contribution to the sum only if its blocks
link the arguments alternatingly in $a$ and $a^{\ast}$.
\\
Case ($1^{\circ}$):
Without loss of generality, we consider the cumulant
\linebreak
$k_n(\dots,b_s,b_{s+1},\dots)$ with $b_s =b_{s+1}=(a^3)^{\ast}$ for
some $s$ with $1 \leq s \leq n-1$. This means that we have to look at
$k_n(\dots, a^{\ast} a^{\ast} a^{\ast}, a^{\ast} a^{\ast}
a^{\ast}, \dots)$. Theorem 2.2 yields in this case
\begin{multline*}
 k_n(\dots, a^{\ast} a^{\ast} a^{\ast},
               a^{\ast} a^{\ast} a^{\ast}, \dots)
   \\ = \sum_{\pi \in NC(3n) \atop \pi \vee \sigma = 1_{3n}}
         k_{\pi} [\dots,a^{\ast},a^{\ast},a^{\ast},
                         a^{\ast},a^{\ast},a^{\ast},\dots] \: ,
\end{multline*}
where $\sigma:=\{\dots,(a^{\ast},a^{\ast},a^{\ast}),
(a^{\ast},a^{\ast},a^{\ast}),\dots\}$. In order to find out
which partitions $\pi \in NC(3n)$ contribute to the sum we look
at the structure of the block containing the element
$b_{s+1,1}=a^{\ast}$; in the following we will call this block
$V$.\\
There are two situations which can occur. The first
possibility is that $b_{s+1,1}$ is the first component of $V$;
in this case the last component of $V$ must be an $a$ and, since
each block has to contain the same number of $a$ and $a^*$, this $a$ has
to be the third $a$ of an argument $a^3$. But then the block $V$ gets
in $\pi\vee\sigma$ not connected with the block containing
$b_{s,3}$ and hence the requirement $\pi\vee\sigma=1_{3n}$ cannot
be fulfilled in such a situation.
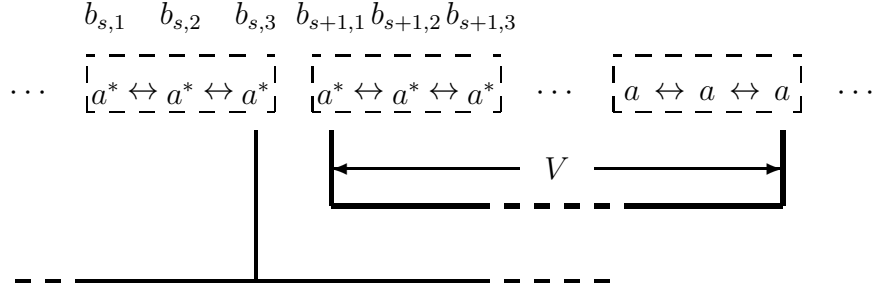
\setlength{\unitlength}{1.0cm}
\[ \begin{picture}(13,4)
   \thicklines
   \put(4,0){\line(0,1){2}}
   \put(1.6,0){\line(1,0){2.4}}
   \put(1.2,0){\line(1,0){0.2}}
   \put(0.8,0){\line(1,0){0.2}}
   \put(4.0,0){\line(1,0){3.1}}
   \put(7.3,0){\line(1,0){0.2}}
   \put(7.7,0){\line(1,0){0.2}}
   \put(8.1,0){\line(1,0){0.2}}
   \put(8.5,0){\line(1,0){0.2}}
   \linethickness{0.6mm}
   \put(5,1){\line(0,1){1}}
   \put(11,1){\line(0,1){1}}
   \put(5.0,1){\line(1,0){2.1}}
   \put(7.3,1){\line(1,0){0.2}}
   \put(7.7,1){\line(1,0){0.2}}
   \put(8.1,1){\line(1,0){0.2}}
   \put(8.5,1){\line(1,0){0.2}}
   \put(8.9,1){\line(1,0){2.1}}
   \put(1,2.5){\makebox(0,0){$\cdots$}}
   \put(2,2.5){\makebox(0,0){$a^{\ast}$}}
   \put(2.5,2.5){\makebox(0,0){$\leftrightarrow$}}
   \put(3,2.5){\makebox(0,0){$a^{\ast}$}}
   \put(3.5,2.5){\makebox(0,0){$\leftrightarrow$}}
   \put(4,2.5){\makebox(0,0){$a^{\ast}$}}
   \put(5,2.5){\makebox(0,0){$a^{\ast}$}}
   \put(5.5,2.5){\makebox(0,0){$\leftrightarrow$}}
   \put(6,2.5){\makebox(0,0){$a^{\ast}$}}
   \put(6.5,2.5){\makebox(0,0){$\leftrightarrow$}}
   \put(7,2.5){\makebox(0,0){$a^{\ast}$}}
   \put(8,2.5){\makebox(0,0){$\cdots$}}
   \put(9,2.5){\makebox(0,0){$a$}}
   \put(9.5,2.5){\makebox(0,0){$\leftrightarrow$}}
   \put(10,2.5){\makebox(0,0){$a$}}
   \put(10.5,2.5){\makebox(0,0){$\leftrightarrow$}}
   \put(11,2.5){\makebox(0,0){$a$}}
   \put(12,2.5){\makebox(0,0){$\cdots$}}
   \thinlines
   \put(1.75,2.25){\dashbox{0.2}(2.5,0.75){}}
   \put(4.75,2.25){\dashbox{0.2}(2.5,0.75){}}
   \put(8.75,2.25){\dashbox{0.2}(2.5,0.75){}}
   \put(2,3.5){\makebox(0,0){$b_{s,1}$}}
   \put(3,3.5){\makebox(0,0){$b_{s,2}$}}
   \put(4,3.5){\makebox(0,0){$b_{s,3}$}}
   \put(5,3.5){\makebox(0,0){$b_{s+1,1}$}}
   \put(6,3.5){\makebox(0,0){$b_{s+1,2}$}}
   \put(7,3.5){\makebox(0,0){$b_{s+1,3}$}}
   \thicklines
   \put(8,1.5){\makebox(0,0){$V$}}
   \put(7.5,1.5){\vector(-1,0){2.5}}
   \put(8.5,1.5){\vector(1,0){2.5}}
   \end{picture}
\]
The second situation that might happen is that $b_{s+1,1}$ is not
the first component of $V$. Then the preceeding element in this block
must be an $a$ and again it must be the third $a$ of an argument $a^3$. But
then the block containing $b_{s,3}$ is again not connected
with $V$ in $\pi\vee\sigma$.
This possibility can be illustrated
as follows:
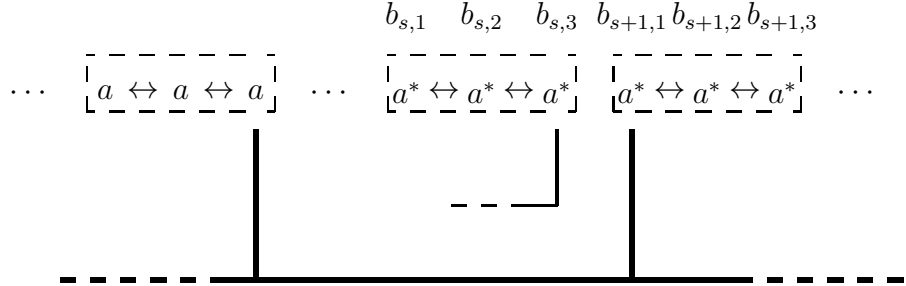
\setlength{\unitlength}{1.0cm}
\[ \begin{picture}(13,4)
\linethickness{0.6mm}
   \put(4,0){\line(0,1){2}}
   \put(1.4,0){\line(1,0){0.2}}
   \put(1.8,0){\line(1,0){0.2}}
   \put(2.2,0){\line(1,0){0.2}}
   \put(2.6,0){\line(1,0){0.2}}
   \put(3.0,0){\line(1,0){0.2}}
   \put(3.4,0){\line(1,0){0.6}}
   \put(4.0,0){\line(1,0){5}}
   \put(9.0,0){\line(0,1){2}}
   \put(9.0,0){\line(1,0){1.6}}
   \put(10.8,0){\line(1,0){0.2}}
   \put(11.2,0){\line(1,0){0.2}}
   \put(11.6,0){\line(1,0){0.2}}
   \put(12.0,0){\line(1,0){0.2}}
   \put(12.4,0){\line(1,0){0.2}}
\thicklines
   \put(6.6,1){\line(1,0){0.2}}
   \put(7.0,1){\line(1,0){0.2}}
   \put(7.4,1){\line(1,0){0.6}}
   \put(8,1){\line(0,1){1}}
   \put(1,2.5){\makebox(0,0){$\cdots$}}
   \put(2,2.5){\makebox(0,0){$a$}}
   \put(2.5,2.5){\makebox(0,0){$\leftrightarrow$}}
   \put(3,2.5){\makebox(0,0){$a$}}
   \put(3.5,2.5){\makebox(0,0){$\leftrightarrow$}}
   \put(4,2.5){\makebox(0,0){$a$}}
   \put(5,2.5){\makebox(0,0){$\cdots$}}
   \put(6,2.5){\makebox(0,0){$a^{\ast}$}}
   \put(6.5,2.5){\makebox(0,0){$\leftrightarrow$}}
   \put(7,2.5){\makebox(0,0){$a^{\ast}$}}
   \put(7.5,2.5){\makebox(0,0){$\leftrightarrow$}}
   \put(8,2.5){\makebox(0,0){$a^{\ast}$}}
   \put(9,2.5){\makebox(0,0){$a^{\ast}$}}
   \put(9.5,2.5){\makebox(0,0){$\leftrightarrow$}}
   \put(10,2.5){\makebox(0,0){$a^{\ast}$}}
   \put(10.5,2.5){\makebox(0,0){$\leftrightarrow$}}
   \put(11,2.5){\makebox(0,0){$a^{\ast}$}}
   \put(12,2.5){\makebox(0,0){$\cdots$}}
   \thinlines
   \put(1.75,2.25){\dashbox{0.2}(2.5,0.75){}}
   \put(5.75,2.25){\dashbox{0.2}(2.5,0.75){}}
   \put(8.75,2.25){\dashbox{0.2}(2.5,0.75){}}
   \put(6,3.5){\makebox(0,0){$b_{s,1}$}}
   \put(7,3.5){\makebox(0,0){$b_{s,2}$}}
   \put(8,3.5){\makebox(0,0){$b_{s,3}$}}
   \put(9,3.5){\makebox(0,0){$b_{s+1,1}$}}
   \put(10,3.5){\makebox(0,0){$b_{s+1,2}$}}
   \put(11,3.5){\makebox(0,0){$b_{s+1,3}$}}
   \end{picture}
\]
Thus, in any case there exists no $\pi$ which fulfills the requirement
$\pi\vee\sigma=1_{3n}$ and hence
$k_n(\dots,a^*a^*a^*,a^*a^*a^*,\dots)$ vanishes in this case.

Case ($2^{\circ}$):
In the case $n$ odd, the cumulant
\linebreak
$k_{\pi}
[b_{1,1},b_{1,2},b_{1,3}, \dots, b_{n,1},b_{n,2},b_{n,3}]$ has a
different number of $a$ and $a^*$ as arguments and hence
at least one of the blocks of $\pi$ cannot be alternating in $a$ and $a^*$.
Thus $k_\pi$ vanishes by the $R$-diagonality of $a$.

As in both cases we do not find any partition giving
a non-vanishing contribution,
the sum vanishes and so do the cumulants $k_n(b_1,\dots,b_n)$.\\
\end{proof}

\subsection{Remark.}
Of course we are now left with the problem of describing the alternating
cumulants of $a^r$ in terms of the alternating cumulants of $a$. We will
provide the solution to this question by showing that the similarity
between $a_1\cdots a_r$ and $a^r$ goes even further as in the Remark 3.10.
Namely, we will show that $a^r$ has the same $*$-distribution
as $a_1\cdots a_r$ if all $a_i$ ($i=1,\dots,r$)
have the same $*$-distribution as $a$. The distribution
of $a^r$ can then be
calculated by an iteration of Prop. 3.9.
In the case of a trace this reduces to a result of Haagerup and
Larsen \cite{HL,L}.
The specical case of powers of a circular element was treated by
Oravecz \cite{O}.

\subsection{Proposition.}
Let $a$ be an $R$-diagonal element and $r$ a positive integer.
Then the $*$-distribution of $a^r$ is the same as the $*$-distribution
of $a_1\cdots a_r$ where each $a_i$ ($i=1,\dots,r$) has the same
$*$-distribution as $a$ and where $a_1,\dots,a_r$ are $*$-free.

\begin{proof}
Since we know that both $a^r$ and $a_1\cdots a_r$ are $R$-diagonal we
only have to see that the respective alternating cumulants coincide.
By Theorem 2.2, we have
\begin{multline*}
k_{2n}(a^r,a^{*r},\dots,a^r,a^{*r})\\=
\sum_{\pi\in NC(2nr)\atop \pi\vee\sigma=1_{2nr}}k_\pi[a,\dots,a,a^*,\dots,
a^*,\dots,a,\dots,a,a^*,\dots,a^*]
\end{multline*}
and
\begin{multline*}
k_{2n}(a_1\cdots a_r,a_r^*\cdots a_1^*,\dots,
a_1\cdots a_r,a_r^*\cdots a_1^*)\\=
\sum_{\pi\in NC(2nr)\atop \pi\vee\sigma=1_{2nr}}
k_\pi[a_1,\dots, a_r,a_r^*,\dots,a_1^*,\dots,
a_1,\dots, a_r,a_r^*,\dots,a_1^*],
\end{multline*}
where in both cases
$\sigma=\{(1,\dots,r),(r+1,\dots,2r),\dots,(2(n-1)r+1,\dots,2nr)\}$.
The only difference between both cases is that in the second case we
also have to take care of the freeness between the $a_i$ which implies
that only such $\pi$ contribute which do not connect different $a_i$.
But the $R$-diagonality of $a$ implies that also in the first case only
such $\pi$ give a non-vanishing contribution, i.e. the freeness in the
second case does not really give an extra condition. Thus both formulas
give the same and the two distributions coincide.
\end{proof}

\end{document}